\documentclass{article}
\usepackage[utf8]{inputenc}
\usepackage{geometry}
\usepackage{amsmath}
\usepackage{amssymb}
\usepackage{amsfonts}
\usepackage{amsthm}
\usepackage{hyperref}
\usepackage{cleveref}
\usepackage{latexsym}
\usepackage{array}
\usepackage{graphicx}
\usepackage{caption}
\usepackage{subcaption}
\usepackage{float}
\usepackage{verbatim}
\usepackage{xcolor}

\newtheorem{theorem}{Theorem}
\newtheorem{definition}{Definition}
\newtheorem{lemma}{Lemma}
\newtheorem{remark}{Remark}
\newtheorem{proposition}{Proposition}

\makeatletter
\usepackage{authblk}
\usepackage{lipsum}
\usepackage{algpseudocode} 

\newcommand{\norm}[1]{\left\lVert#1\right\rVert}
\newcommand{\brac}[1]{\left(#1\right)}
\newcommand{\avg}[1]{\left<#1\right>}

\title{Uniformly higher order accurate schemes for dynamics of charged particles under fast oscillating magnetic fields\footnote{Submitted to arXiv on November 7, 2024. \\ \textbf{Funding:} The first author was funded by the Ministry of Education, Government of India under the Prime Minister's Research Fellowship/grant PM/MHRD-19-17567.03. The last author benefits from fundings by the European Union via the Euratom Research and Training Programme (Grant Agreement No 101052200 EUROFusion). Views and opinions expressed are however those of the author only and do not necessarily reflect those of the European Union or the European Commission. Neither the European Union nor the European Commission can be held responsible for them.}}
\author[1]{\small \textbf{Megala Anandan}.} 
\author[2]{\small \textbf{Benjamin Boutin}}
\author[3]{\small \textbf{Nicolas Crouseilles}} 
\affil[1]{Indian Institute of Science, C.V. Raman Road, 560012, Bangalore, India. mail: megalaa@iisc.ac.in.}
\affil[2]{Univ Rennes, CNRS, IRMAR - UMR CNRS 6625, F-35000 Rennes, France. mail: benjamin.boutin@univ-rennes.fr.}
\affil[3]{Universit\'e de Rennes, Inria Rennes (Mingus team) and IRMAR UMR CNRS 6625, F-35042 Rennes, France \& ENS Rennes. mail: nicolas.crouseilles@inria.fr.}
\date{}

\begin{document}

\maketitle
\begin{abstract}

This work deals with the numerical approximation of   plasmas which are confined by the effect of a fast oscillating magnetic field (see~\cite{Bostan2012}) in the Vlasov model.  
The presence of this magnetic field induces  oscillations (in time) to the solution of the characteristic equations. 
Due to its multiscale character, a standard time discretization would lead to an inefficient solver. In this work, time integrators are derived and analyzed 
for a class of highly oscillatory differential systems. We prove the uniform accuracy property of these time integrators, meaning that the accuracy does not depend on the small parameter $\varepsilon$. Moreover, we construct an extension of the scheme which degenerates towards an energy preserving numerical scheme for the averaged model, when $\varepsilon\to 0$.  
Several numerical results illustrate the capabilities of the method.
\end{abstract}

\section{Introduction}
The confinement of charged particles system is of important interest nowadays, in particular 
due to the application to fusion devices. The main strategy used in the machines like tokamaks or stellarators is to apply a strong magnetic field ${\bf B}(t, x)$ ($t$ denotes the time variable and $x$ the space variable) to confine the particles far from the boundary wall of the device. Indeed, due to the extreme temperature of plasma, it is crucial to avoid direct contact 
of the particles with wall material while allowing 
the plasma to reach fusion conditions. 
To confine particles, it is also possible to consider fast oscillating (in time) magnetic fields instead of strong ones. Indeed, the choice proposed in~\cite{Bostan2012} is the following: ${\bf B}(t, x) = \theta(t/\varepsilon) B(x)b(x), 0<\varepsilon <\!\!< 1$, 
where $\theta$ is a given $P$-periodic function, 
$B(x)$ is a scalar positive function, and $b(x)$ is a unitary vector. This configuration may require less energy to design compared to the 
strong magnetic field configuration and it has been proven that this configuration has good confinement properties  (see~\cite{Bostan2012}). 

In this work, we are interested in the numerical approximation of the model obtained with this highly-oscillatory magnetic field. From the Vlasov equation~\eqref{vlasov} in a four-dimensional setting in phase space, the external magnetic field 
is supposed to be ${\bf B}(t, x) = \theta(t/\varepsilon) (0, 0, B)^T$ with $B>0$. The electric field has to be decomposed into a gradient and rotational part so that we can formulate the corresponding characteristic equations as follows: 
\begin{eqnarray}
\label{carac_intro_x}
\dot{x}(t) &=& v(t), \\ 
\label{carac_intro_v}
\dot{v} (t) &=& E(t, x(t)) +\textcolor{black}{\frac{B}{2 \varepsilon} \theta'\brac{t/\varepsilon}Jx(t)}  + \theta(t/\varepsilon) B J v(t), \;\; t\in [0, T]
\end{eqnarray}
where $x(t)\in\mathbb{R}^2$ denotes the position, $v(t)\in\mathbb{R}^2$ denotes the velocity of the particle, $E:\mathbb{R}_+\times \mathbb{R}^2\to \mathbb{R}^2$ denotes a given irrotational  electric field (deriving from an electric potential $E=-\nabla \phi$, with $\phi:\mathbb{R}_+\times \mathbb{R}^2\to \mathbb{R}$) 
and $J$ denotes the canonical $2\times 2$ symplectic matrix. 

Our goal is to construct and analyse efficient numerical schemes to approximate differential systems of the form \eqref{carac_intro_x}-\eqref{carac_intro_v}. 
It is worth noticing that the system~\eqref{carac_intro_x}-\eqref{carac_intro_v} 
is highly oscillatory in time ; in particular the second order derivatives of $v$ is stiff which 
makes the derivation of higher order schemes challenging 
since the consistency error usually involves the time derivatives of the solution. To capture the correct dynamics with standard numerical strategies, the time step $\Delta t$ has to be smaller than $\varepsilon\in (0,1]$, leading to very expensive methods. 
Thus, the first goal is to overcome the stiffness induced by the highly oscillatory character of the solution. In addition, the method has to be consistent with the asymptotic (or averaged) model 
obtained in the limit $\varepsilon\to 0$ (see ~\cite{Bostan2012}). 
These properties are usually referred in the literature as asymptotic-preserving or uniformly accurate strategies which ensures that the limit of the scheme is a scheme for the limit model. 

During the last decade, several works have been proposed to solve the highly-oscillatory 
equations efficiently and accurately. These schemes have been referred as uniformly accurate since the error can be proved to be uniform with respect to $\varepsilon$. One can quote the Duhamel based strategies applied for Klein-Gordon or Dirac equations \cite{schratz_ua, ua_cai, ua_ode}, two-scale methods used for Vlasov equations \cite{two_scale_cpcsever, two_scale_flr, two_scale_fsl, two_scale_jcp1,  two_scale_jcp2, two_scale_jcp3d, two_scale_jin, two_scale_mathcomp, two_scale_hoode, two_scale_numer, two_scale_prel, ua_geom}, micro-macro approach \cite{ua_gilles,two_scale_jin,ua_derivativefree}, 
integrator scheme for advection-diffusion equations \cite{astuto, astuto2}, multiscale decomposition or exponential wave integrators for Klein-Gordon, Dirac or Schrodinger equations  \cite{bao_ua_ode1,bao_ua_ode2,bao_ua_ode3,bao_ua_ode4,bao_ua_ode5}. 
In this work, we also consider a Duhamel based approach to design uniformly accurate numerical scheme for~\eqref{carac_intro_x}-\eqref{carac_intro_v}. 
A first order uniformly accurate scheme can be easily  obtained by freezing the unknown in time and integrating the highly oscillating part exactly, which mainly consists of integrating the function $\theta$. 
To derive higher order numerical methods, it is possible to recursively substitute adequately the Duhamel formula in itself, in the spirit of \cite{ua_cai, astuto}. 

In the linear context, we will consider the following general framework 
\begin{equation}
\label{linear_intro}
\dot{U}(t) = A(t/\varepsilon)U(t)
\end{equation}
with $U(t)\in \mathbb{R}^d$ and $A$ a given $d\times d$ matrix-valued $P$-periodic function. In this case, 
it is possible 
to derive arbitrarily higher order uniformly accurate explicit schemes. On top of that, a suitable use of the Duhamel formula enables to design (second order) midpoint uniformly accurate numerical schemes. Such schemes are particularly relevant for the configuration where the system \eqref{linear_intro}  has a specific structure (Hamiltonian or $A$ skew-symmetric). Indeed, in this case, the averaged model derived from~\eqref{linear_intro} becomes $\dot{\bar{U}}(t) = \langle A\rangle \bar{U}(t)$ when $\varepsilon\to 0$, and midpoint method 
is the method of choice since it preserves some invariants ($L^2$ norm or energy) exactly. Hence, our goal is to propose numerical schemes that  are not only uniformly accurate, but also degenerate towards a midpoint scheme for the averaged model when $\varepsilon\to 0$ .

Then, we investigate how to generalize the strategy to nonlinear systems of the form 
$$
\dot{U}(t) = A(t/\varepsilon) U(t) + g(t, U(t)).
$$
Even if the first order case can be easily obtained, the extension to higher order requires more attention. 
Indeed, substituting the Duhamel formula in the nonlinear term leads to a highly oscillatory nonlinear term in $g$ which is not tractable from a numerical point of view. A suitable expansion allows the oscillatory term to be explicitly revealed, enabling the development of higher order, uniformly accurate explicit schemes in the nonlinear context. However, the extension to midpoint strategy (which usually offers invariant preserving properties) turns out to be too costly in view of Vlasov applications since it involves nonlinear fixed point techniques. To overcome this drawback, we employ the SAV technique introduced in~\cite{sav2} to construct efficient structure-preserving numerical methods for nonlinear problems. After reformulating the system~\eqref{carac_intro_x}-\eqref{carac_intro_v} using the SAV framework, we propose  an efficient uniformly accurate midpoint numerical scheme 
which degenerates (when $\varepsilon\to 0$) 
towards a midpoint scheme for the averaged model 
of~\eqref{carac_intro_x}-\eqref{carac_intro_v}. 
It turns out that the averaged model enjoys a Hamiltonian structure (see~\cite{Bostan2012}). Hence, a midpoint (with an explicit complexity) numerical scheme is particularly relevant since it exactly preserves the discrete energy. 


The rest of the paper is organized as follows. 
In Section \ref{models}, the models we are interested in are presented, together with their asymptotic limit. 
In section \ref{Section:UA numerical scheme}, the numerical schemes are presented in a  linear 
context: higher order explicit schemes and second order midpoint schemes are introduced and their uniform accuracy is proved. In section \ref{Section:Energy preserving scheme for avg model}, we extend the numerical scheme to the nonlinear case. In Section \ref{pic}, we explain how these solvers can be integrated in Particle-In-Cell framework to approximate the nonlinear Vlasov-Poisson system. 
Finally,  several numerical results are shown to illustrate the properties of the schemes in different contexts. 

\section{Models}
\label{models}

The main motivation comes from a physical application 
introduced in~\cite{Bostan2012} where the author studies 
the effect of a fast oscillating external magnetic field 
on a system of charged particles. Let us thus consider the Vlasov equation with the fast oscillating magnetic field ${\bf B}(t) = \theta(t/\varepsilon)(0, 0, B)^T$ (with $B>0$). A more general model would include gradient and curvature effects, which will be studied in a future work. But since the magnetic field depends on time, the electric field induced from the Maxwell equations takes the form $\frac{B}{2\varepsilon} \theta'(t/\varepsilon) Jx$ (see~\cite{Bostan2012}). 
Hence, the Vlasov equation becomes  
\begin{gather}
\label{vlasov}
    \partial_t f^{\varepsilon} + v \cdot \nabla_x f^{\varepsilon} + \brac{ E(t,x) + \frac{B}{2 \varepsilon} \theta'\brac{t/\varepsilon}Jx + B \theta\brac{t/\varepsilon}Jv } \cdot \nabla_v f^{\varepsilon} = 0, \\
    f^{\varepsilon} (0,x,v) = f^{\text{in}} (x,v), \;\; (x,v) \in \mathbb{R}^2 \times \mathbb{R}^2, \;\; t\in [0, T], 
\end{gather}
where $f^{\varepsilon}(t,x,v)\geq 0$ is the distribution function of particles in the phase space and $E=-\nabla \phi$ is a given irrotational electric field. Here, $B$ is a constant, $\theta(s)$ is a given $P$-periodic function of class ${\cal C}^1$, and $0<\varepsilon \ll 1$. Further, for any $a=(a_1,a_2)^T \in \mathbb{R}^2$, $Ja = J \begin{bmatrix} a_1, a_2 \end{bmatrix}^T = \begin{bmatrix} a_2, -a_1 \end{bmatrix}^T$. \\
As a transport equation, one has $\frac{d}{dt} f(t, x(t), v(t)) = 0$ which translates the fact $f$ is constant along the characteristic curves. Considering $B=1$, $(x(t), v(t))$ is solution to the following ODE system:
\begin{equation}
\label{ODE(x,v)}
    \begin{gathered}
    \dot{x} = v, \\
    \dot{v} = E(t,x) + \frac{1}{2 \varepsilon} \theta'\brac{t/\varepsilon}Jx +  \theta\brac{t/\varepsilon}Jv.
    \end{gathered}
\end{equation} 
This system is not only stiff but also involves highly oscillatory 
terms making its numerical approximation quite challenging 
since a restrictive condition on the time step is commonly required. Before 
discussing the construction of appropriate numerical schemes for system~\eqref{ODE(x,v)}, 
let us study its asymptotic behavior. As suggested in~\cite{Bostan2012}, we first perform the 
following change of variable $q:= v - \frac{1}{2} \theta\brac{t/\varepsilon}Jx$, so as to reformulate~\eqref{ODE(x,v)} as a system in the unknown $(x,q)$:
\begin{equation}
\label{ODE(x,q)_original}
    \begin{gathered}
    \dot{x} = q +  \frac{1}{2} \theta\brac{t/\varepsilon}Jx, \\
    \dot{q} = E(t,x) + \frac{1}{2} \theta\brac{t/\varepsilon}Jq + \frac{1}{4}\theta^2\brac{t/\varepsilon}J^2 x.
\end{gathered}
\end{equation}
Thanks to this change of variable, the system has 
the form of a standard highly-oscillatory problem 
for which the limit $\varepsilon\to 0$ is well known~\cite{allaire, Bogolyubov1961AsymptoticMI}. Namely, introducing the averaged value $\avg{\theta}$ of $\theta$  (with $\avg{\theta}:=\frac{1}{P}\int_0^P \theta (s) ds$), the solution $(x,q)$ of~\eqref{ODE(x,q)_original} converges when $\varepsilon\to 0$ towards the solution $(\bar{x}, \bar{q})$ of the limiting system:
\begin{equation}
\label{averaged_ODE}
    \begin{gathered}
    \dot{\bar{x}} = \bar{q} +  \frac{1}{2} \langle\theta\rangle J\bar{x}, \\
    \dot{\bar{q}} = E(t,\bar{x}) + \frac{1}{2}  \langle\theta\rangle J\bar{q} + \frac{1}{4}  \langle\theta^2\rangle J^2 \bar{x}. 
\end{gathered}
\end{equation}

As observed in~\cite{Bostan2012}, this averaged model~\eqref{averaged_ODE} enjoys a Hamiltonian structure where  the Hamiltonian is 
$H(x,q)=\frac{1}{2}|q|^2 +  \frac{1}{2}\langle\theta\rangle q\cdot Jx + \frac{1}{8}\langle\theta^2\rangle |x|^2 + \phi(x)$ so that the system~\eqref{averaged_ODE} can be reformulated as  
$$
  \begin{gathered}
\left(
\begin{matrix}
    \dot{\bar{x}}\\
    \dot{\bar{q}} 
    \end{matrix} 
  \right) 
  = 
  \left(
\begin{matrix}
   0 & I \\
-I & 0 
    \end{matrix} 
  \right) 
\left(
\begin{matrix}
    \nabla_x H(\bar{x}, \bar{q})\\
    \nabla_q H(\bar{x}, \bar{q})
        \end{matrix} 
  \right),  \;\;\;\;\;\;  
  I=\left(\begin{matrix}
    1 & 0\\
    0 & 1
        \end{matrix} 
  \right). 
\end{gathered}
$$ 

In the rest of the paper, for the sake of simplicity, we will consider the following normalized model which shares the same properties as the 'physical' one 
\eqref{ODE(x,q)_original} 
\begin{equation}
\label{ODE(x,q)}
    \begin{gathered}
    \dot{x} = q +  \theta\brac{t/\varepsilon}Jx, \\
    \dot{q} = E(t,x) + \theta\brac{t/\varepsilon}Jq +  \theta^2\brac{t/\varepsilon}J^2 x.
\end{gathered}
\end{equation}
For the model \eqref{ODE(x,q)}, the 
averaged model is characterized by the following (normalized) Hamiltonian which writes 
\begin{equation}
\label{ham_norm}
H(\bar{x}, \bar{q}) = \frac{1}{2}|\bar{q}|^2 +  \langle\theta\rangle \bar{q}\cdot J\bar{x} + \frac{1}{2}\langle\theta^2\rangle |\bar{x}|^2 + \phi(\bar{x}). 
\end{equation}

In this work, we are interested in the construction of efficient numerical schemes for more generel highly oscillatory models of the form
\begin{equation}
\label{lin_general}
\dot{U} = A(t/\varepsilon) U + g(U), \;\;\; U(t=0)=U_0\in \mathbb{R}^d,  
\end{equation}
where  $A$ is $d\times d$ matrix-valued $P$-periodic function and $g:\mathbb{R}^d \mapsto \mathbb{R}^d$ is a smooth given function. It is known (see~\cite{allaire, Bogolyubov1961AsymptoticMI}) that the solution $U(t)$ of~\eqref{lin_general} converges towards 
$\bar{U}(t)$ which is the solution of 
\begin{equation}
\label{lin_general_limit}
\dot{\bar{U}} = \langle A\rangle {\bar{U}} + g(\bar{U}), \;\;\; \mbox { where } \langle A\rangle = \frac{1}{P}\int_0^P A(s) ds \;\; \mbox{ and } \;\;  \bar{U}(t=0)=\bar{U}_0.  
\end{equation}
As an example, the system~\eqref{ODE(x,q)} involves the following oscillating $4\times 4$ matrix valued function $A(t/\varepsilon)$ and its average $\langle A\rangle$ defined by  
\begin{equation}
    A(t/\varepsilon) = \begin{pmatrix}
        \theta(t/\varepsilon)J & I\\
        -\theta^2(t/\varepsilon)I & \theta(t/\varepsilon)J
    \end{pmatrix} \;\;\; \mbox{ and } \;\;\; 
    \langle A \rangle = \begin{pmatrix}
 \langle \theta \rangle J & I\\
        - \langle \theta^2 \rangle I & \langle \theta \rangle J
    \end{pmatrix}. 
\end{equation}
Let us remark that in the linear case $g=0$, the averaged model \eqref{lin_general_limit} reduces to a simple linear system $\dot{\bar{U}} = \langle A\rangle \bar{U}$, whose dynamics is determined by the eigenvalues of $\langle A\rangle$. This will be useful to compare with the numerical results.

\section{Uniformly accurate numerical schemes: linear case}
\label{Section:UA numerical scheme}
In this section, we present numerical discretizations of 
highly oscillatory linear ODEs of the form 
\begin{equation}
\label{ODEsys_E0}
\dot{U}(t) = A(t/\varepsilon) U(t),  \;\;\; U(t=0)=U_0,  
\end{equation}
where $U: \mathbb{R}_+\mapsto \mathbb{R}^d$, $U_0\in \mathbb{R}^d$ and $A: \mathbb{R}_+\mapsto {\mathcal M}_{d,d}(\mathbb{R})$ is $P$-periodic. Let us recall that when $\varepsilon$ tends to $0$, the solution to~\eqref{ODEsys_E0} converges towards the solution to the following averaged model 
\begin{equation}
\label{ODEsys_E0_averaged}
\dot{\bar{U}}(t) = \langle A\rangle \bar{U}(t),  \;\;\; \bar{U}(t=0)=U_0. 
\end{equation}
For the time discretization, we will consider the time step 
$\Delta t>0$ and the following notations for $n\in\mathbb{N}$: $t_n=n\Delta t$ and $U_n\in\mathbb{R}^d$ that denotes the numerical solution at time $t_n$, hopefully an approximation of $U(t_n)$. 

\subsection{Arbitrarily higher order accurate scheme: explicit case} 
\label{lin_explicit}

Before presenting the general framework, let us introduce the Duhamel formula for the solution of~\eqref{ODEsys_E0} step by step, by integrating it between $t_n$ and $t_{n+1}$: 
\begin{equation}
\label{duhamel_1}
U(t_{n+1}) = U(t_n) + \int_{t_n}^{t_{n+1}} A\brac{s_1/\varepsilon} U(s_1) ds_1. 
\end{equation} 
Integrating~\eqref{ODEsys_E0} again, but 
between $t_n$ and $s_1$ gives: 
$$
U(s_1) = U(t_n) + \int_{t_n}^{s_1} A\brac{s_2/\varepsilon} U(s_2) ds_2, 
$$
and substituting this in~\eqref{duhamel_1} leads to 
\begin{eqnarray*}
U(t_{n+1}) &=& U(t_n) + \int_{t_n}^{t_{n+1}} A\brac{s_1/\varepsilon} \Big[U(t_n) + \int_{t_n}^{s_1} A\brac{s_2/\varepsilon} U(s_2) ds_2 \Big] ds_1 \nonumber\\
& =& U(t_n) + \int_{t_n}^{t_{n+1}} A\brac{s_1/\varepsilon} ds_1 U(t_n) + 
 \int_{t_n}^{t_{n+1}} A\brac{s_1/\varepsilon} \int_{t_n}^{s_1} A\brac{s_2/\varepsilon} U(s_2) ds_2 ds_1. \nonumber
\end{eqnarray*}
Let $p\geq 1$ be a fixed integer, it is of course possible to recursively substitute $(p-1)$ times the Duhamel formula into itself to get  
\begin{eqnarray}
    U(t_{n+1}) 
    &=& U(t_n) + \int_{t_n}^{t_{n+1}} A\brac{s_1/\varepsilon} ds_1 U(t_n) + \int_{t_n}^{t_{n+1}} A\brac{s_1/\varepsilon} \int_{t_n}^{s_1} A\brac{s_2/\varepsilon} U(s_2) ds_2 ds_1 \nonumber\\
    &=& U(t_n) + \int_{t_n}^{t_{n+1}} A\brac{s_1/\varepsilon} ds_1 U(t_n) + \int_{t_n}^{t_{n+1}} A\brac{s_1/\varepsilon} \int_{t_n}^{s_1} A\brac{s_2/\varepsilon}  ds_2 ds_1 U(t_n) \nonumber\\ && + \int_{t_n}^{t_{n+1}} A\brac{s_1/\varepsilon} \int_{t_n}^{s_1} A\brac{s_2/\varepsilon} \int_{t_n}^{s_2} A\brac{s_3/\varepsilon} U(s_3) ds_3 ds_2 ds_1 \nonumber\\
    &=& U(t_n) + \sum_{k=1}^{p-1} \int_{t_n}^{t_{n+1}} \int_{t_n}^{s_1} \dots \int_{t_n}^{s_{k-1}} A\brac{s_1/\varepsilon} A\brac{s_2/\varepsilon} \dots A\brac{s_k/\varepsilon} ds_k ds_{k-1}\dots ds_1 \ U(t_n) \nonumber\\
    \label{Arbit order exact ODE_(x,q)}
    && + \int_{t_n}^{t_{n+1}} \int_{t_n}^{s_1} \dots \int_{t_n}^{s_{p-1}} A\brac{s_1/\varepsilon} A\brac{s_2/\varepsilon} \dots A\brac{s_p/\varepsilon} U(s_p) ds_p ds_{p-1}\dots ds_1.
\end{eqnarray}
Thus, 
the nested Duhamed formula~\eqref{Arbit order exact ODE_(x,q)} suggests the 
following $p$-order numerical scheme  for $U_n\approx U(t_n)$: 
\begin{align}
\hspace{-0.25cm}U_{n+1}&\hspace{-0.1cm}= U_n + \sum_{k=1}^{p-1} \int_{t_n}^{t_{n+1}} \int_{t_n}^{s_1} \dots \int_{t_n}^{s_{k-1}} A\brac{s_1/\varepsilon} A\brac{s_2/\varepsilon} \dots A\brac{s_k/\varepsilon} ds_k ds_{k-1}\dots ds_1 \ U_n  \nonumber\\
    \label{Arbit order ODE_(x,q)}
    & + \int_{t_n}^{t_{n+1}} \int_{t_n}^{s_1} \dots \int_{t_n}^{s_{p-1}} A\brac{s_1/\varepsilon} A\brac{s_2/\varepsilon} \dots A\brac{s_p/\varepsilon} ds_p ds_{p-1}\dots ds_1 \ U_n, \\
     &\hspace{-0.25cm}= \! U_n \! + \Big(\sum_{k=1}^{p} H_k \Big) U_n, \;\; \mbox{with } \;\;  H_k =\!\! \int_{t_n}^{t_{n+1}}\!\!\!\! \int_{t_n}^{s_1} \!\!\dots \! \int_{t_n}^{s_{k-1}} \!\!\!\! A\brac{s_1/\varepsilon} A\brac{s_2/\varepsilon} \dots A\brac{s_k/\varepsilon} ds_k ds_{k-1}\dots ds_1. 
\end{align}
\begin{remark}
    When required, the integrals (with respect to $t$) of products of the periodic function $A\brac{t/\varepsilon}$, that are involved in the matrices $H_m$, are computed analytically. 
    When difficult or even impossible, the function can be expressed through Fourier series $\brac{\theta(t/\varepsilon)}^m = \sum_{k \in \mathbb{Z}} \hat{\theta}_{m_k} e^{ikt/\varepsilon}$ and suitable quadrature $\brac{\theta(t/\varepsilon)}^m \simeq \sum_{|k| \leq K} \hat{\theta}_{m_k} e^{ikt/\varepsilon}$ can be used for their approximations. 
\end{remark}

The numerical scheme~\eqref{Arbit order ODE_(x,q)} is uniformly accurate of order $p$, as stated in the following theorem.  
\begin{theorem}
%
    Let $U(t=0) \in \mathbb{R}^d$. and $A\brac{s=t/\varepsilon} \in {\mathcal M}_{d,d}(\mathbb{R})$ be a given bounded $P$-periodic function of $s$. Consider the solution $t\mapsto U(t)\in\mathbb{R}^d$  to~\eqref{ODEsys_E0} over a given time interval $[0,T]$. The numerical scheme~\eqref{Arbit order ODE_(x,q)} with initial data $U_0=U(t=0)$ is uniformly accurate of order $p$:\\
    for $\Delta t>0$, denote $N$ the largest integer such that $t_N\leq T$, then 
    $$\norm{U(t_n)-U_n} \leq C \Delta t^p$$ for all $n=1,2,\dots,N$, with a constant $C$ independent of $\Delta t$ and $\varepsilon$.
\end{theorem}
\begin{proof}
From Lemma \ref{appendix_tech_lemma}-\ref{appendix_Gronwall}, $U(t)$ is bounded over the time interval $[0,T]$, uniformly in $\varepsilon$. 
 We consider the difference between~\eqref{Arbit order exact ODE_(x,q)} and~\eqref{Arbit order ODE_(x,q)} to obtain an induction on the error, as in \cite{astuto} 
    \begin{eqnarray*}
        e_{n+1} &=& e_n + \sum_{m=1}^{p-1} H_m e_n + \int_{t_n}^{t_{n+1}} \int_{t_n}^{s_1} \dots \int_{t_n}^{s_{p-1}} A\brac{s_1/\varepsilon} A\brac{s_2/\varepsilon} \dots A\brac{s_p/\varepsilon} \brac{U(s_p)-U_n} ds_p ds_{p-1}\dots ds_1 \\
        &=& \brac{I+\sum_{m=1}^{p-1} H_m + H_p} e_n \\ && + \int_{t_n}^{t_{n+1}} \int_{t_n}^{s_1} \dots \int_{t_n}^{s_{p-1}} A\brac{s_1/\varepsilon} A\brac{s_2/\varepsilon} \dots A\brac{s_p/\varepsilon} \brac{U(s_p)-U(t_n)} ds_p ds_{p-1}\dots ds_1 \\
        &=& \brac{I+\sum_{m=1}^{p} H_m } e_n \\
        && + \int_{t_n}^{t_{n+1}} \int_{t_n}^{s_1} \dots \int_{t_n}^{s_{p-1}} \int_{t_n}^{s_p} A\brac{s_1/\varepsilon} A\brac{s_2/\varepsilon} \dots A\brac{s_p/\varepsilon}A\brac{s_{p+1}/\varepsilon} U(s_{p+1}) ds_{p+1}ds_p ds_{p-1}\dots ds_1
    \end{eqnarray*}
    where $t_n\leq s_{p+1} \leq s_p$. Further since $A\brac{t/\varepsilon}$ and ${U(t)}$ are uniformly bounded in $\varepsilon$, we may denote the finite quantities $C_m = \norm{A}_{L^\infty}^m$ for $1\leq m\leq p$, so as to bound
the norm of $e_{n+1}$  by $\norm{e_{n+1}} \leq \brac{1+\sum_{m=1}^p C_m \Delta t^m} \norm{e_n} + C \Delta t^{p+1}$, 
    and we conclude by the Gronwall lemma. 
\end{proof}

In addition to the previous theorem, the following proposition identifies the asymptotic numerical scheme, which is obtained by taking the limit in~\eqref{Arbit order ODE_(x,q)}. As expected, thanks to the UA property, this asymptotic scheme is a $p$-th order accurate scheme for~\eqref{ODEsys_E0_averaged}. 
\begin{proposition}
The numerical scheme~\eqref{Arbit order ODE_(x,q)} converges when $\varepsilon\to 0$ to the $p$-th order accurate scheme below, consistent with the limiting equation~\eqref{ODEsys_E0_averaged}:
$$
{U}_{n+1} =\Big( \sum_{k=0}^{p} \frac{\Delta t^k}{k !} \langle A\rangle ^k \Big) {U}_{n}.
$$
\end{proposition}
\begin{proof}
Using lemma \ref{appendix_tech_lemma}-\ref{lemma_limit}, we first get 
$\int_{t_n}^{s_{m-1}} A(t/\varepsilon) dt = (s_{m-1}-t_n) \langle A \rangle + \mathcal{O}({\varepsilon})$  
and then
\begin{eqnarray*}
\int_{t_n}^{s_{m-2}} \int_{t_n}^{s_{m-1}} A(s_{m-1}/\varepsilon) A(s_m/\varepsilon)   ds_{m}ds_{m-1}  &=& \int_{t_n}^{s_{m-2}} A(s_{m-1}/\varepsilon) \int_{t_n}^{s_{m-1}}A(s_m/\varepsilon)   ds_{m}ds_{m-1} \nonumber\\
&\hspace{-12cm}= & \hspace{-6cm}
\int_{t_n}^{s_{m-2}} A(s_{m-1}/\varepsilon) (s_{m-1} - t_n) ds_{m-1} \langle A \rangle + \mathcal{O}({\varepsilon}) = \frac{(s_{m-2}-t_n)^2}{2} \langle A \rangle^2 + \mathcal{O}({\varepsilon}).  
\end{eqnarray*}
By induction, we similarly obtain
    \begin{eqnarray*}
        H_m &=& \int_{t_n}^{t_{n+1}} \int_{t_n}^{s_1} \dots \int_{t_n}^{s_{m-1}} A\brac{s_1/\varepsilon} A\brac{s_2/\varepsilon} \dots A\brac{s_m/\varepsilon} ds_m ds_{m-1}\dots ds_1 \nonumber\\
        &=& \int_{t_n}^{t_{n+1}} \frac{(s_1-t_n)^ {m-1}}{(m-1)!} ds_1 \langle A \rangle^m + {\cal O}(\varepsilon) =  \frac{\Delta t^m}{m!} \langle A \rangle^m + {\cal O}(\varepsilon). 
    \end{eqnarray*} 
Therefore, the scheme~\eqref{Arbit order ODE_(x,q)} becomes when $\varepsilon\to 0$: ${U}_{n+1} = {U}_{n} +  \Big( \sum_{k=1}^{p} \frac{\Delta t^k}{k!} \langle A\rangle ^k \Big) {U}_{n} + {\cal O}(\varepsilon)$. 
Passing to the limit $\varepsilon\to 0$ leads to the result. 
\end{proof}

\subsection{Uniformly accurate numerical schemes: midpoint case}
In this part, we investigate midpoint schemes combined with the strategy presented in the previous subsection. 
The reason for studying a possible extension of explicit schemes 
to midpoint schemes lies in the structure of the asymptotic model. 
Indeed, under some assumptions on the matrix $\langle A \rangle$  (skew symmetry or hamiltonian structure), the asymptotic model preserves some invariants ($L^2$ norm or some energies) 
and midpoint schemes are known to preserve the discrete counterpart. 
Hence, our goal is to construct a numerical scheme that 
is uniformly accurate and degenerates towards a scheme (for the averaged model) that also preserves the discrete invariants, as $\varepsilon\to 0$. 


Before introducing our strategy, let us first discuss 
a direct (but naive) way to construct midpoint schemes 
for~\eqref{ODEsys_E0}. 
A direct extension of the explicit scheme presented in the 
previous subsection would be 
\begin{equation}
\label{mid_point_naive}
U_{n+1} = U_n+ \int^{t_{n+1}}_{t_n} A(s/\varepsilon) ds 
\frac{U_{n+1} + U_n}{2}. 
\end{equation}
However, even if such a scheme converges towards a midpoint scheme 
for the averaged model~\eqref{ODEsys_E0_averaged}, it is only first  order uniformly accurate  (even if it is second order accurate for any fixed $\varepsilon$). 
See Appendix \ref{appendix1}. 

We now improve the previous scheme~\eqref{mid_point_naive} so as to get a second order UA scheme that degenerates towards a midpoint scheme for the averaged model~\eqref{averaged_ODE}.  
On one side, we have the forward Duhamel formula 
\begin{equation}
\label{duhamel_exp}
U(s)=U(t_n) + \int_{t_n}^s A(s_1/\varepsilon) U(s_1) ds_1, \mbox{ for } s\geq t_n, 
\end{equation}
and on the other side, we have the backward Duhamel formula 
\begin{equation}
\label{duhamel_imp}
U(s)=U(t_{n+1}) - \int^{t_{n+1}}_s A(s_1/\varepsilon) U(s_1) ds_1, \mbox{ for } s\leq t_{n+1}.  
\end{equation}
Looking back at the derivation of the UA scheme of previous subsection, we write 
\begin{equation}
\label{duhamel_imex}
U(t_{n+1}) = U(t_n)  + \int^{t_{n+1}}_{t_n} A(s/\varepsilon) U(s) ds,
\end{equation}
and the strategy was to replace $U(s)$ by the corresponding Duhamel formula 
which enabled us to get explicit UA schemes. To derive midpoint UA schemes, 
we now instead replace $U(s)$ by the average of~\eqref{duhamel_exp} and~\eqref{duhamel_imp}: 
\begin{equation}
\label{us_midpoint}
U(s) = \frac{U(t_{n+1}) +U(t_n)}{2} + \frac{1}{2}\int_{t_n}^s A(s_1/\varepsilon) U(s_1) ds_1 -\frac{1}{2}\int^{t_{n+1}}_s A(s_1/\varepsilon) U(s_1) ds_1, 
\end{equation}
so that~\eqref{duhamel_imex} becomes 
\begin{eqnarray}
U(t_{n+1}) &=& U(t_n)  + \int^{t_{n+1}}_{t_n} A(s/\varepsilon) ds \frac{U(t_{n+1}) +U(t_n)}{2} \nonumber\\
\label{exact_midpoint}
&&\hspace{-1cm}+ \frac{1}{2} \int^{t_{n+1}}_{t_n} A(s/\varepsilon)\int_{t_n}^s A(s_1/\varepsilon) U(s_1) ds_1 ds -\frac{1}{2}\int^{t_{n+1}}_{t_n} A(s/\varepsilon)
\int^{t_{n+1}}_s A(s_1/\varepsilon) U(s_1) ds_1 ds. 
\end{eqnarray}
Considering the numerical unknown $U_n\approx U(t_n)$, one has the following second order UA numerical scheme 
\begin{eqnarray}
U_{n+1} &=& U_n + \int^{t_{n+1}}_{t_n} A(s/\varepsilon) ds \frac{U_{n+1} +U_n}{2}\nonumber\\
\label{midpoint}
&& \hspace{-1cm}+ \frac{1}{2} \Big(\int^{t_{n+1}}_{t_n} A(s/\varepsilon)\int_{t_n}^s A(s_1/\varepsilon)  ds_1 ds- \int^{t_{n+1}}_{t_n} A(s/\varepsilon)
\int^{t_{n+1}}_s A(s_1/\varepsilon) ds_1 ds\Big)\frac{U_{n+1} +U_n}{2}. 
\end{eqnarray}

\begin{theorem}
    Let $U(t=0) \in \mathbb{R}^d$ and $A\brac{s=t/\varepsilon} \in {\mathcal M}_{d,d}(\mathbb{R})$ be a given bounded $P$-periodic function of $s$. Consider the solution $t\mapsto U(t)\in\mathbb{R}^d$  to~\eqref{ODEsys_E0} over a given time interval $[0,T]$. The numerical scheme~\eqref{midpoint} with initial data $U_0=U(t=0)$ is second order uniformly accurate:\\
    for $\Delta t>0$, denote $N$ the largest integer such that $t_N\leq T$, then 
    $$\norm{U(t_n)-U_n} \leq C \Delta t^2$$ for all $n=1,2,\dots,N$, with a constant $C$ independent of $\Delta t$ and of $\varepsilon$.
\end{theorem}
\begin{proof}
Here again, from Lemma \ref{appendix_tech_lemma}-\ref{appendix_Gronwall}, $U(t)$ is bounded over the time interval $[0,T]$, uniformly in $\varepsilon$.\\
Let us find an induction formula on the error $e_{n+1} = U(t_{n+1}) - U_{n+1}$ by considering the difference between~\eqref{exact_midpoint} and~\eqref{midpoint} 
\begin{eqnarray*}
e_{n+1} &=& e_n  + \frac{1}{2}\int^{t_{n+1}}_{t_n} A(s/\varepsilon)ds \Big[ e_{n+1} + e_n \Big] +  \frac{1}{2}\int^{t_{n+1}}_{t_n} A(s/\varepsilon)\int_{t_n}^s A(s_1/\varepsilon) \Big[U(s_1) - \frac{U(t_{n+1}) +U(t_n)}{2}   \Big] ds_1 ds \nonumber\\
&& - \frac{1}{2}\int^{t_{n+1}}_{t_n} A(s/\varepsilon)\int_s^{t_{n+1}} A(s_1/\varepsilon) \Big[U(s_1) - \frac{U(t_{n+1}) +U(t_n)}{2}   \Big] ds_1 ds\nonumber\\ 
 &&+\frac{1}{2}\int^{t_{n+1}}_{t_n} A(s/\varepsilon)\int_{t_n}^s A(s_1/\varepsilon) \Big[\frac{e_{n+1} +e_n}{2}   \Big] ds_1 ds  - \frac{1}{2}\int^{t_{n+1}}_{t_n} A(s/\varepsilon)\int_s^{t_{n+1}} A(s_1/\varepsilon) \Big[ \frac{e_{n+1} +e_n}{2}   \Big] ds_1 ds. 
\end{eqnarray*}
From~\eqref{us_midpoint}, the third term (and similarly for the fourth term) becomes 
\begin{eqnarray*}
 &&\frac{1}{2}\int^{t_{n+1}}_{t_n} A(s/\varepsilon)\int_{t_n}^s A(s_1/\varepsilon) \Big[U(s_1) - \frac{U(t_{n+1}) +U(t_n)}{2}   \Big] ds_1 ds \nonumber\\
 &&=  \frac{1}{4}\int^{t_{n+1}}_{t_n} A(s/\varepsilon)\int_{t_n}^s A(s_1/\varepsilon) 
 \Big[\int_{t_n}^{s_1} A(s_2/\varepsilon) U(s_2)ds_2 - \int_{s_1}^{t_{n+1}}  A(s_2/\varepsilon) U(s_2)ds_2   \Big]ds_1 ds 
\end{eqnarray*} 
so that the error induction becomes 
\begin{eqnarray*}
e_{n+1} &=& e_n  \\ 
&&\hspace{-1.4cm} + \Big(\frac{1}{2}\int^{t_{n+1}}_{t_n} \!\! A(s/\varepsilon)ds  +\frac{1}{4}\int^{t_{n+1}}_{t_n} \!\!A(s/\varepsilon)\int_{t_n}^s A(s_1/\varepsilon) ds_1 ds  - \frac{1}{4} \int^{t_{n+1}}_{t_n} \!\!A(s/\varepsilon)\int_s^{t_{n+1}} \!\!A(s_1/\varepsilon)  ds_1 ds \Big)(e_{n+1} +e_n) \\
&&\hspace{-1.4cm} +  \frac{1}{4}\int^{t_{n+1}}_{t_n} A(s/\varepsilon)\int_{t_n}^s A(s_1/\varepsilon) 
 \Big[\int_{t_n}^{s_1} A(s_2/\varepsilon) U(s_2)ds_2 - \int_{s_1}^{t_{n+1}}  A(s_2/\varepsilon) U(s_2)ds_2   \Big]ds_1 ds  \\
 &&\hspace{-1.4cm} - \frac{1}{4}\int^{t_{n+1}}_{t_n} A(s/\varepsilon)\int_{s}^{t_{n+1}} A(s_1/\varepsilon) 
 \Big[\int_{t_n}^{s_1} A(s_2/\varepsilon) U(s_2)ds_2 - \int_{s_1}^{t_{n+1}}  A(s_2/\varepsilon) U(s_2)ds_2   \Big]ds_1 ds.  
\end{eqnarray*}
From the uniform boundedness of $U$, the terms in the two last lines are uniformly bounded by 
$C\Delta t^3$ since the intervals of integration are all smaller than $\Delta t$. One can then write 
$$
\|e_{n+1}\| \leq \|e_n\| + (C\Delta t + C\Delta t^2) (\|e_{n+1}\| + \|e_{n}\|) + C\Delta t^3,  
$$
or $\brac{1-C\Delta t} \|e_{n+1}\| \leq \brac{1+C\Delta t} \|e_{n}\| + C\Delta t^3$.  
Thus, we obtain (for $\Delta t < 1/C$) 
$$
\|e_{n+1}\| \leq \frac{1+C\Delta t}{1-C\Delta t} \|e_{n}\| + \frac{C}{1-C\Delta t} \Delta t^3 \leq  (1+K\Delta t) \|e_{n}\| + C \left(1+K\Delta t\right) \Delta t^3 \leq (1+K\Delta t) \|e_{n}\| + K \Delta t^3, 
$$
and we use the Gronwall lemma to conclude that the error is uniformly 
bounded by $C\Delta t^2$. We used the simple estimates 
$(1+C\Delta t)/(1-C\Delta t) \leq 1+K\Delta t$ and $1/(1-C\Delta t) \leq 1+K\Delta t$. 
\end{proof}

As done in the previous subsection, we now identify 
the asymptotic numerical scheme, that is the numerical
scheme obtained by passing to the limit in~\eqref{midpoint}. 

\begin{proposition}
The midpoint scheme~\eqref{midpoint} converges when $\varepsilon\to 0$ to the following numerical scheme 
\begin{equation}
\label{midpoint_limit}
\bar{U}_{n+1} = \bar{U}_n + \Delta t \langle A\rangle  \frac{\bar{U}_n + \bar{U}_{n+1}}{2}, 
\end{equation}
which is a second order approximation of the solution $\bar{U}(t_{n+1}) = \exp(\langle A\rangle (n+1) \Delta t) U_0$ of~\eqref{ODEsys_E0_averaged}. 
Moreover, when $\langle A\rangle =-\langle A\rangle^T$, the scheme~\eqref{midpoint_limit} preserves the discrete $L^2$ norm $\|\bar{U}_{n}\| = \|\bar{U}_0\|$ for all $n\geq 0$. Finally, for $\bar{U}_0-U_0= O (\varepsilon)$ and any time $T>0$, there exists a constant $C>0$ independent of $\varepsilon$ and $\Delta t$ such that provided $0\leq n\Delta t\leq T$, one has
\[ \|\bar{U}_n-U_n\| \leq C \varepsilon.\]
\end{proposition}
\begin{proof}
By using Lemma~\ref{appendix_tech_lemma}-\ref{lemma_limit}, the scheme~\eqref{midpoint} is interpreted as
\begin{eqnarray*}
{U}_{n+1} &=& {U}_n \!+ \Delta t \langle A\rangle \frac{{U}_n + {U}_{n+1}}{2} \!+\! \left(\frac{\Delta t^2}{4} \langle A\rangle^2 \!-\! \frac{\Delta t^2}{4} \langle A\rangle^2 \right) \frac{U_n + {U}_{n+1}}{2} \!+ \mathcal{O}(\Delta t\,\varepsilon)\\
&=& {U}_n \!+\! \Delta t \langle A\rangle \frac{{U}_n + {U}_{n+1}}{2} \!+ \mathcal{O}(\Delta t\,\varepsilon),  
\end{eqnarray*}
thereby proving the first part of the proposition. For the second part, we consider the scalar product of~\eqref{midpoint_limit} with $(\bar{U}_{n+1} + \bar{U}_n)$ and use the skew-symmetry of $\langle A\rangle$.  
The last point of the proof is obtained by comparing the schemes~\eqref{midpoint} and~\eqref{midpoint_limit}. Setting
$$C_n(\varepsilon)= \frac{1}{2} \Big(\int^{t_{n+1}}_{t_n} A(s/\varepsilon)\int_{t_n}^s A(s_1/\varepsilon)  ds_1 ds- \int^{t_{n+1}}_{t_n} A(s/\varepsilon)
\int^{t_{n+1}}_s A(s_1/\varepsilon) ds_1 ds\Big),$$
and observing then from Lemma~\ref{appendix_tech_lemma}-\ref{lemma_limit}, that $C_n(\varepsilon) = \mathcal{O}(\Delta t\,\varepsilon)$, we have 
\begin{eqnarray*}
\bar{U}_{n+1}-U_{n+1} &=& (\bar{U}_{n}-U_{n}) + \Delta t \langle A\rangle \bar{U}_{n+1/2} - \int_{t_n}^{t_{n+1}}A(s/\varepsilon)ds U_{n+1/2} - C_n(\varepsilon) U_{n+1/2}\\
&=& (\bar{U}_{n}-U_{n}) + \Delta t \langle A\rangle (\bar{U}_{n+1/2}-U_{n+1/2}) + (\Delta t \langle A\rangle - \textstyle\int_{t_n}^{t_{n+1}}A(s/\varepsilon)ds) U_{n+1/2} - C_n(\varepsilon) U_{n+1/2}\\
&=& (I-\Delta t \langle A\rangle /2)^{-1}\Big((I+\Delta t \langle A\rangle /2)(\bar{U}_{n}-U_{n}) + \mathcal{O}(\Delta t\,\varepsilon) U_{n+1/2}\Big).
\end{eqnarray*}
From the boundedness of the numerical solution in bounded times, uniform in $\varepsilon$, and Gronwall inequalities, we then prove the result.
\end{proof}
\begin{remark}
The scheme \eqref{midpoint} enjoys strong links with Magnus methods \cite{iserles}. Indeed, 
in our case, the second order Magnus method would write 
$U_{n+1} =\exp(\Theta) U_n$, 
with $\Theta=\int_{t_n}^{t_{n+1}}A(t/\varepsilon)dt + \frac{1}{2}\int_{t_n}^{t_{n+1}}\int_{t_n}^t [A(t), A(s)]dsdt$ (with the commutator bracket $[A(t), A(s)] = A(t)A(s)-A(s)A(t)$, which vanishes for scalar $A$ $(d=1)$). A 
second order expansion of this scheme gives $\exp(\Theta) \approx I+\Theta+\frac{1}{2}\Theta^2$ which corresponds a second order expansion of \eqref{midpoint}. 
This link motivates the use of Lie group methods \cite{iserles} in our context due to their geometric properties in more general applications. 
As an example, 
the Magnus scheme  
$U_{n+1}= \exp(\Theta) U_n$ can be proved to be second order UA, and if $A$ is skew-symmetric, it also preserves the $L^2$ norm of the numerical solution. This class of schemes will be studied in a future work. 
\end{remark}

\subsection{Midpoint schemes for the charged particles model: linear case}
In this part, we focus on the charged particles model by applying 
the scheme~\eqref{midpoint} to the model~\eqref{ODE(x,q)_original}. 
Then, we obtain 
\begin{eqnarray}
\label{midpoint_particles_eq}
\left(
\begin{matrix} 
x_{n+1}\\
q_{n+1} 
\end{matrix}
\right) 
&=& 
\left(
\begin{matrix} 
x_{n}\\
q_{n} 
\end{matrix}
\right)
+
{\cal B}
\left(
\begin{matrix} 
\frac{x_{n+1} + x_n}{2}\\
\frac{q_{n+1} + q_n}{2} 
\end{matrix}
\right)+ \frac{1}{2}{\cal A}\left(
\begin{matrix} 
\frac{x_{n+1} + x_n}{2}\\
\frac{q_{n+1} + q_n}{2} 
\end{matrix}
\right),  
\end{eqnarray}
where ${\cal B}$ and ${\cal A}$ are  defined by 
\begin{equation}
\label{midpoint_B}
{\cal B}=\left(
\begin{matrix} 
\int_{t_n}^{t_{n+1}} \theta(s/\varepsilon) ds J & \Delta t I\\
\int_{t_n}^{t_{n+1}} \theta^2(s/\varepsilon) ds J^2 & \int_{t_n}^{t_{n+1}} \theta(s/\varepsilon) ds J
\end{matrix}
\right), \;\;\;\; {\cal A}=\left(
\begin{matrix} 
{\cal A}_{1,1}&{\cal A}_{1,2}\\
{\cal A}_{2,1}&{\cal A}_{2,2}
\end{matrix}
\right)\in {\mathcal M}_{4,4}(\mathbb{R}),
\end{equation}
where the block components ${\cal A}_{i,j}\in {\mathcal M}_{2,2}(\mathbb{R})$ are given by 
\begin{eqnarray}
\label{a11}
\hspace{-5.5cm}{\cal A}_{1,1} \!\!\!\!\! &=&\!\!\!\!\!  
\int_{t_n}^{t_{n+1}}\!\!\!\! \theta(s/\varepsilon)\Big[\!\!-\!\!\int_{t_n}^{s}  \!\!\theta(s_1/\varepsilon) ds_1  \!+\! \int^{t_{n+1}}_{s}  \!\!\!\!\theta(s_1/\varepsilon) ds_1 \Big] ds I  \!-\! \int_{t_n}^{t_{n+1}} \!\!\Big[\int_{t_n}^{s}  \!\!\theta^2(s_1/\varepsilon) ds_1  \!-\!  \int^{t_{n+1}}_{s} \!\!\!\! \theta^2(s_1/\varepsilon) ds_1 \Big] ds I \\ 
\label{a12}
{\cal A}_{1,2} \!\!\!\!&\!\!\!\!=\!\!\!\!&\!\!\!\! 
\int_{t_n}^{t_{n+1}} \!\!\!\!\theta(s/\varepsilon)((s-t_n)-(t_{n+1}-s_1)) ds J + \int_{t_n}^{t_{n+1}} \Big[\int_{t_n}^{s}  \theta(s_1/\varepsilon) ds_1 -\int^{t_{n+1}}_{s}  \theta(s_1/\varepsilon) ds_1 \Big] ds J \\
\label{a21}
{\cal A}_{2,1} \!\!\!\!\!\!&=&\!\!\!\!\!\! 
\int_{t_n}^{t_{n+1}} \!\!\!\!\theta^2(s/\varepsilon) \Big[\!\!-\!\!\int_{t_n}^{s}  \!\!\theta(s_1/\varepsilon) ds_1  \!+\! \int^{t_{n+1}}_{s}  \!\!\!\!\theta(s_1/\varepsilon) ds_1    \Big] ds J  \!\!-\!\! \int_{t_n}^{t_{n+1}}\!\!\!\!\theta(s/\varepsilon)\Big[\int_{t_n}^{s}\!\!  \theta^2(s_1/\varepsilon) ds_1 \! -\!  \int^{t_{n+1}}_{s}  \!\!\!\!\theta^2(s_1/\varepsilon) ds_1 \Big] ds J \\
\label{a22}
{\cal A}_{2,2} \!\!\!\!&=&\!\!\!\! -\int_{t_n}^{t_{n+1}}\!\!\theta^2(s/\varepsilon)((s-t_n)\!-\!(t_{n+1}-s_1)) ds I  + \int_{t_n}^{t_{n+1}}\!\!\!\!\theta(s/\varepsilon)
\Big[\!-\!\int_{t_n}^{s}  \theta(s_1/\varepsilon) ds_1  + \int^{t_{n+1}}_{s} \!\!\!\! \theta(s_1/\varepsilon) ds_1    \Big] ds I . 
\end{eqnarray}
In the following theorem, we identify the asymptotic numerical scheme by taking the limit $\varepsilon\to 0$ in~\eqref{midpoint_particles_eq}. 
The asymptotic model is satisfied by $(\bar{x}, \bar{q})$ 
solution to 
\begin{eqnarray}
\label{averaged_sec31}
    \dot{\bar{x}} &=& \bar{q} + \langle \theta\rangle J\bar{x}, \\
    \label{averaged_sec32}
    \dot{\bar{q}} &=&  \langle \theta\rangle J\bar{q}+\langle \theta^2\rangle J^2\bar{x}.   
\end{eqnarray}
This system has a Hamiltonian structure 
with a Hamiltonian given by \eqref{ham_norm}, 
and this can be written as $H(x,q) = H_1(x,q)+H_2(x,q)$ defined by 
\begin{equation} 
\label{def_subh}
H_1(x, q) = \frac{1}{2}|q|^2 + \frac{1}{2}\langle\theta^2\rangle |x|^2, \;\mbox{ and } \; H_2(x, q) = \langle\theta\rangle q\cdot Jx. 
\end{equation} 
Moreover, we have $\frac{d}{dt}H_1(x(t),q(t)) :=\frac{d}{dt} \Big(\frac{1}{2}|q(t)|^2 + \frac{1}{2}\langle\theta^2\rangle |x(t)|^2\Big) = \frac{d}{dt}H_2(x(t),q(t))=\frac{d}{dt} \Big( \langle\theta\rangle q(t)\cdot Jx(t)\Big) = 0$ where $(x,q)$ is solution to \eqref{averaged_sec31}-\eqref{averaged_sec32}.

\begin{theorem}
\label{theorem_h1h2}
The scheme~\eqref{midpoint_particles_eq} is second order UA and degenerates  when $\varepsilon\to 0$ towards an energy preserving scheme for the averaged model \eqref{averaged_sec31}-\eqref{averaged_sec32}, whose solution $(\bar{x}_n, \bar{q}_n)_n$ satisfies: $H_1(\bar{x}_{n+1}, \bar{q}_{n+1})=H_1(\bar{x}_{n}, \bar{q}_{n})$ 
and $H_2(\bar{x}_{n+1}, \bar{q}_{n+1})=H_2(\bar{x}_{n}, \bar{q}_{n})$ where $H_1$ and $H_2$ 
are given by~\eqref{def_subh}. 
\end{theorem}

\begin{proof}
The UA property follows from the previous theorem. In this proof, we focus on the asymptotic behavior of the 
scheme~\eqref{midpoint_particles_eq} when $\varepsilon\to 0$. 
We use the estimates from Lemma~\ref{appendix_tech_lemma}-\ref{lemma_limit} several times, with a crude upper bound of the form $O(\varepsilon)$ (instead of $O(\varepsilon\Delta t)$ or $O(\varepsilon^2)$) so as to highlight the formal limiting scheme. When $\varepsilon\to 0$, the matrix ${\cal B}$ satisfies  
$$
{\cal B} = \left(
\begin{matrix} 
\int_{t_n}^{t_{n+1}} \theta(s/\varepsilon) ds J & \Delta t I\\
\int_{t_n}^{t_{n+1}} \theta^2(s/\varepsilon) ds J^2 & \int_{t_n}^{t_{n+1}} \theta(s/\varepsilon) ds J
\end{matrix}
\right) 
= 
\left(
\begin{matrix}
\Delta t \langle \theta \rangle J & \Delta t I\\
\Delta t \langle \theta^2 \rangle J^2 & \Delta t \langle \theta \rangle J 
\end{matrix}\right) + \mathcal{O}(\varepsilon). 
$$
For the matrix ${\cal A}$, let us consider each block component. 
For the first one, we get from Lemma~\ref{appendix_tech_lemma}-\ref{lemma_limit} 
\begin{eqnarray*}
{\cal A}_{1,1} &=&
\Big( 
\int_{t_n}^{t_{n+1}}\theta(s/\varepsilon) \int_{t_n}^s \theta(s_1/\varepsilon)ds_1 ds - \int_{t_n}^{t_{n+1}}\theta(s/\varepsilon) \int_s^{t_{n+1}} \theta(s_1/\varepsilon)ds_1 ds\nonumber\\
&&+ \int_{t_n}^{t_{n+1}} \int_{t_n}^s \theta^2(s_1/\varepsilon)ds_1 ds 
- \int_{t_n}^{t_{n+1}} \int_s^{t_{n+1}} \theta^2(s_1/\varepsilon)ds_1 ds \Big) J^2\nonumber\\
&=&  \frac{\Delta t^2}{2} \langle \theta\rangle^2 -\frac{\Delta t^2}{2} \langle \theta\rangle^2  +\frac{\Delta t^2}{2} \langle \theta^2\rangle -\frac{\Delta t^2}{2} \langle \theta^2\rangle +{\cal O}(\varepsilon)={\cal O}(\varepsilon). 
\end{eqnarray*}
For the second component, we get from Lemma~\ref{appendix_tech_lemma}-\ref{lemma_limit} 
\begin{eqnarray*}
{\cal A}_{1,2} &=& \Big( 
\int_{t_n}^{t_{n+1}} \theta(s/\varepsilon) ((s-t_n)-(t_{n+1}-s)) ds + \int_{t_n}^{t_{n+1}} \int_{t_n}^s \theta(s_1/\varepsilon)ds_1 ds  - 
\int_{t_n}^{t_{n+1}} \int_s^{t_{n+1}} \theta(s_1/\varepsilon)ds_1 ds\Big) J\nonumber\\
&=&  
\frac{\Delta t^2}{2}\langle\theta\rangle - \frac{\Delta t^2}{2}\langle\theta\rangle + \frac{\Delta t^2}{2}\langle\theta\rangle -\frac{\Delta t^2}{2}\langle\theta\rangle +{\cal O}(\varepsilon)= {\cal O}(\varepsilon).  
\end{eqnarray*}
And similarly, we have for the component ${\cal A}_{2,1}={\cal O}(\varepsilon)$ and 
$ {\cal A}_{2,2}={\cal O}(\varepsilon)$.  

Hence, the numerical scheme obtained when $\varepsilon \to 0$ (with $x_n\to \bar{x}_n$ and $q_n \to \bar{q}_n$) becomes 
\begin{eqnarray}
\label{sec3:midpoint_linear}
\left(
\begin{matrix} 
\bar{x}_{n+1}\\
\bar{q}_{n+1} 
\end{matrix}
\right) 
&=& 
\left(
\begin{matrix} 
\bar{x}_{n}\\
\bar{q}_{n} 
\end{matrix}
\right) + 
\left(
\begin{matrix} 
\Delta t \langle \theta \rangle J & \Delta t I\\
\Delta t \langle \theta^2 \rangle J^2 & \Delta t \langle \theta \rangle J
\end{matrix}
\right)
\left(
\begin{matrix} 
\frac{\bar{x}_{n+1} + \bar{x}_n}{2}\\
\frac{\bar{q}_{n+1} + \bar{q}_n}{2} 
\end{matrix}
\right). 
\end{eqnarray}
It remains to prove that the scheme~\eqref{sec3:midpoint_linear} preserves the energies $H_1$ and $H_2$. 

Let us multiply the first equation of~\eqref{sec3:midpoint_linear} by $\bar{x}_{n+1/2}=(\bar{x}_{n+1}+\bar{x}_n)/2$ to get 
$$
(|\bar{x}_{n+1}|^2 - |\bar{x}_{n}|^2)/\Delta t = \bar{x}_{n+1/2} \cdot \bar{q}_{n+1/2},  
$$
since $x\cdot Jx=0$. Further, multiplying the second equation by $\bar{q}_{n+1/2}$ gives 
$$
(|\bar{q}_{n+1}|^2 - |\bar{q}_{n}|^2)/\Delta t =  - \langle\theta^2\rangle  \bar{q}_{n+1/2} \cdot \bar{x}_{n+1/2}. 
$$
Combining the two last equalities leads to 
$|\bar{q}_{n+1}|^2 + \langle\theta^2\rangle |\bar{x}_{n+1}|^2 = |\bar{q}_{n}|^2 +  \langle\theta^2\rangle |\bar{x}_{n}|^2$ which corresponds to the $H_1$ preservation. Moreover, multiplying the first equation of~\eqref{sec3:midpoint_linear} by $J$ and considering the scalar product by $\bar{q}_{n+1/2}$ enables to get 
\begin{eqnarray*}
\bar{q}_{n+1/2}\cdot J  (\bar{x}_{n+1} - \bar{x}_n)/\Delta t  &=&  \langle\theta\rangle \bar{q}_{n+1/2}\cdot J J{\bar{x}}_{n+1/2} + \bar{q}_{n+1/2}\cdot J \bar{q}_{n+1/2}= - \langle\theta\rangle \bar{q}_{n+1/2}\cdot {\bar{x}}_{n+1/2}.
\end{eqnarray*}
Now, multiplying the second equation of~\eqref{sec3:midpoint_linear} by $J \bar{x}_{n+1/2}$ leads to 
$$
J \bar{x}_{n+1/2} \cdot  (\bar{q}_{n+1} - \bar{q}_n)/\Delta t  = \langle\theta\rangle J \bar{x}_{n+1/2} \cdot J \bar{q}_{n+1/2} =\langle\theta\rangle  \bar{x}_{n+1/2} \cdot  \bar{q}_{n+1/2}.   
$$
Adding the two last equalities leads to  
\begin{eqnarray*}
\bar{q}_{n+1/2}\cdot J  (\bar{x}_{n+1} - \bar{x}_n) + J \bar{x}_{n+1/2} \cdot  (\bar{q}_{n+1} - \bar{q}_n)&=&0, 
\end{eqnarray*}
but the left hand side can be reformulated as 
\begin{eqnarray}
\bar{q}_{n+1/2}\cdot J  (\bar{x}_{n+1} - \bar{x}_n) + J \bar{x}_{n+1/2} \cdot  (\bar{q}_{n+1} - \bar{q}_n)&=&
(\bar{q}_{n+1}+\bar{q}_{n})\cdot J  (\bar{x}_{n+1} - \bar{x}_n) + J (\bar{x}_{n+1}+\bar{x}_n) \cdot  (\bar{q}_{n+1} - \bar{q}_n)\nonumber\\
&& \hspace{-7.5cm} =\bar{q}_{n+1} \cdot J  \bar{x}_{n+1} -\bar{q}_{n+1} \cdot J  \bar{x}_{n} +\bar{q}_{n} \cdot J  \bar{x}_{n+1}  - \bar{q}_n\cdot J \bar{x}_n 
+ J \bar{x}_{n+1}\cdot \bar{q}_{n+1} - J\bar{x}_{n+1} \cdot \bar{q}_n + J\bar{x}_n \cdot \bar{q}_{n+1}  - J\bar{x}_n\cdot \bar{q}_n \nonumber\\
\label{relation_discrete}
&&\hspace{-7.5cm}  = \bar{q}_{n+1} \cdot J  \bar{x}_{n+1} - \bar{q}_n\cdot J \bar{x}_n, 
\end{eqnarray}
which corresponds to the $H_2$ preservation and concludes the proof. 
\end{proof}

\begin{remark}
It is actually sufficient to prove that the scheme~\eqref{midpoint_particles_eq} degenerates towards a midpoint scheme for the averaged model \eqref{averaged_sec31}-\eqref{averaged_sec32}, when $\varepsilon\to 0$. Indeed, 
since $H_1$ and $H_2$ are quadratic, we can appeal Theorem 2.1, p101 of~\cite{gni} to ensure that the midpoint method preserves quadratic invariants. However, the detailed proof will be useful for the nonlinear case. 
\end{remark}


\section{Nonlinear system for charged particles under fast oscillating magnetic field}
In this section, we consider some possible extensions of the numerical schemes introduced in previous section to the nonlinear 
case.  

\subsection{Explicit case}
We first extend the explicit uniformly accurate numerical schemes presented in subsection \ref{lin_explicit} to nonlinear case. In this subsection, we consider nonlinear ODEs of the form 
\begin{equation}
\label{ODE(x,q)_vector}
    \dot{U}(t) = A(t/\varepsilon) U(t) + g(t,U(t)), \;\;\; U(t=0)=U_0.  
\end{equation} 
for which we propose first and second order uniformly accurate schemes.  

\subsubsection{First order accurate scheme}
Let us consider the general form for nonlinear systems of the form \eqref{ODE(x,q)_vector} for which the following 
first order uniformly accurate numerical scheme can be obtained from the previous section  
\begin{equation}
\label{First order_ODE(x,q)}
    U_{n+1} = U_n + \int_{t_n}^{t_{n+1}} A\brac{t/\varepsilon } dt\, U_n + \Delta t g_n \text{ with } g_n= g\brac{t_n,U_n}.
\end{equation}
We have the following theorem. 
\begin{theorem}
    Let $U(t=0) \in \mathbb{R}^d$, $A\brac{s=t/\varepsilon} \in {\mathcal M}_{d,d}(\mathbb{R})$ be a given bounded $P$-periodic function of $s$, and $g$ be a given bounded Lipschitz continuous function of $(t,U)$ with Lipschitz constant $K$. 
    Consider the solution $t\mapsto U(t)\in\mathbb{R}^d$  to~\eqref{ODE(x,q)_vector} over a given time interval $[0,T]$.
    The numerical scheme~\eqref{First order_ODE(x,q)} with initial data $U_0=U(t=0)$ is first order uniformly accurate:\\
    for $\Delta t>0$, denote $N$ the largest integer such that $t_N\leq T$, then 
    $$\norm{U(t_n)-U_n} \leq C \Delta t$$ for all $n=1,2,\dots,N$, with a constant $C$ independent of $\Delta t$ and of $\varepsilon$.
\end{theorem}
\begin{proof}
    Integrating~\eqref{ODE(x,q)_vector} between $t_n$ and $t$ gives,
    \begin{equation*}
    \label{First order exact subs_ODE(x,q)}
        U(t)=U(t_n) + \int_{t_n}^{t} A\brac{s/\varepsilon} U(s) ds  + \int_{t_n}^{t} g(s,U(s)) ds. 
    \end{equation*}
    Evaluating the above exact integration at $t=t_{n+1}$, we get
    \begin{equation}
    \label{First order exact_ODE(x,q)}
        U(t_{n+1})=U(t_n) + \int_{t_n}^{t_{n+1}} A\brac{s/\varepsilon} U(s) ds  + \int_{t_n}^{t_{n+1}} g(s,U(s)) ds. 
    \end{equation}
    Then with $e_n=U(t_n)-U_n$, the difference between \eqref{First order_ODE(x,q)} and \eqref{First order exact_ODE(x,q)} gives
    \begin{eqnarray*}
        e_{n+1} &=& e_n + \int_{t_n}^{t_{n+1}} A\brac{s/\varepsilon } \brac{U(s)-U_n} ds + \int_{t_n}^{t_{n+1}} \brac{g(s,U(s))-g_n} ds \\
        &=& e_n + \int_{t_n}^{t_{n+1}} A\brac{s/\varepsilon} \brac{U(s)-U(t_n)} ds + \int_{t_n}^{t_{n+1}} A\brac{s/\varepsilon} e_n ds + \int_{t_n}^{t_{n+1}} \brac{g(s,U(s))-g_n} ds \\
        &=& \brac{I + \int_{t_n}^{t_{n+1}} A\brac{s/\varepsilon} ds} e_n + \int_{t_n}^{t_{n+1}} A\brac{s/\varepsilon} \int_{t_n}^{s} \Big[ A\brac{s_1/\varepsilon} U(s_1)  +  g(s_1,U(s_1)) \Big] ds_1 ds \\ 
        & & + \int_{t_n}^{t_{n+1}} \brac{g(s,U(s))-g_n} ds.
    \end{eqnarray*}
    Further since $\norm{A\brac{t/\varepsilon}}$, $\norm{U(t)}$, and $\norm{g(t,U(t))}$ are uniformly bounded and 
    and we also have $\norm{g(s,U(s))-g_n} \leq K \norm{(s,U(s))-(t_n,U_n)} \leq C (\Delta t + \norm{e_n})$, where $K$ is the Lipschitz constant.
    Then the norm of $e_{n+1}$ becomes
    \begin{eqnarray*}
        \norm{e_{n+1}} &\leq& \brac{1+ C \Delta t } \norm{e_n} + C \Delta t^2 + C \Delta t \brac{ \Delta t  + \norm{e_n}}\leq  \brac{1+C \Delta  t}  \norm{e_n} + C \Delta t^2, 
    \end{eqnarray*}
    where the constant $C$ can changed from one line to the next but it is independent of $\Delta t,\varepsilon$. 
    Thus we conclude by Gronwall lemma that $\norm{e_n} \leq C \Delta t$ with $C$ independent of $\Delta t$ and $\varepsilon$, for $n=1,2,..,N$. 
\end{proof}

\subsubsection{Second order accurate scheme}
In this section, we present how to extend the above first order scheme to  get a second order uniformly accurate scheme for \eqref{ODE(x,q)_vector}. The following second order uniformly accurate scheme 
(and the corresponding theorem) holds for $g(t,U(t))$ but for the sake of simplicity, we will present only for the case $g(U(t))$.  

To derive a second order scheme \eqref{ODE(x,q)_vector}, 
we follow the strategy introduced in the linear case: we integrate \eqref{ODE(x,q)_vector} between $t_n$ and $t_{n+1}$ and perform a recursive substitution to get 
\begin{eqnarray}
    U(t_{n+1})&=& U(t_n) + \int_{t_n}^{t_{n+1}} A\brac{s_1/\varepsilon} U(s_1) ds_1  + \int_{t_n}^{t_{n+1}} g(U(s_1)) ds_1 \nonumber \\
    &=& U(t_n) + \int_{t_n}^{t_{n+1}} A\brac{s_1/\varepsilon} ds_1 U(t_n) + \int_{t_n}^{t_{n+1}} A\brac{s_1/\varepsilon} \int_{t_n}^{s_1}  A\brac{s_2/\varepsilon} U(s_2) ds_2 ds_1 \nonumber \\
    \label{Second order exact_Enon0_ODE_(x,q)}
    && + \int_{t_n}^{t_{n+1}} A\brac{s_1/\varepsilon} \int_{t_n}^{s_1} g(U(s_2)) ds_2 ds_1  + \int_{t_n}^{t_{n+1}} g(U(t_n) + h(s_1))ds_1, 
\end{eqnarray}
with $h(s_1):=\int_{t_n}^{s_1} \Big[ A\brac{s_2/\varepsilon} U(s_2) + g(U(s_2)) \Big] ds_2$. 
In view of the numerical scheme, the solution is frozen at time $t_n$ in the third and fourth terms whereas we have to discuss the last term. First,  
assuming $g$ is smooth, we perform a Taylor expansion of $g(U(t_n) + h(s_1))\approx g(U(t_n))+(\nabla g)(U(t_n)) h(s_1)$ and second, one uses a suitable approximation $\tilde{h}$ of $h$. 
We then obtain the following second order numerical discretization of \eqref{ODE(x,q)_vector}  
\begin{eqnarray}
    U_{n+1} &=& U_n + \int_{t_n}^{t_{n+1}} A\brac{s_1/\varepsilon} ds_1 U_n + \int_{t_n}^{t_{n+1}} A\brac{s_1/\varepsilon} \int_{t_n}^{s_1}  A\brac{s_2/\varepsilon} ds_2 ds_1 U_n \nonumber \\
    \label{Second order_Enon0_ODE_(x,q)}
    &&  + \int_{t_n}^{t_{n+1}} A\brac{s_1/\varepsilon} \int_{t_n}^{s_1}  ds_2 ds_1 g_n + \Delta t g_n + \nabla \left. g \right|_n\int_{t_n}^{t_{n+1}} \tilde{h}(s_1) ds_1  
\end{eqnarray}
where $g_n = g\brac{U_n}$ and $\nabla \left. g \right|_n = (\nabla  g){\brac{U_n}}$ and $\tilde{h}$ is an approximation of $h$ given by 
$$
\tilde{h}(s_1) =  \int_{t_n}^{s_1} A\brac{s_2/\varepsilon} ds_2 U_n + (s_1-t_n) g_n.  
$$  
\begin{theorem}
\label{theorem_order2_nl}
    Let $U(t=0) \in \mathbb{R}^d$, $A\brac{s=t/\varepsilon} \in {\mathcal M}_{d,d}(\mathbb{R})$ be a given bounded $P$-periodic function of $s$, and $g$, $\nabla g$ be two given bounded Lipschitz continuous functions of $U$ with Lipschitz constant $K$. 
    Consider the solution $t\mapsto U(t)\in\mathbb{R}^d$  to~\eqref{ODE(x,q)_vector} over a given time interval $[0,T]$.
    The numerical scheme~\eqref{Second order_Enon0_ODE_(x,q)} with initial data $U_0=U(t=0)$ is second order uniformly accurate:\\
    for $\Delta t>0$, denote $N$ the largest integer such that $t_N\leq T$, then 
    $$\norm{U(t_n)-U_n} \leq C \Delta t^2$$ for all $n=1,2,\dots,N$, with a constant $C$ independent of $\Delta t$ and of $\varepsilon$.
\end{theorem}
\begin{proof}
    Considering the notation $e_n=U(t_n)-U_n$ for $n=1,2,\dots,N$, subtracting \eqref{Second order exact_Enon0_ODE_(x,q)} with \eqref{Second order_Enon0_ODE_(x,q)} gives
    \begin{eqnarray*}
        e_{n+1} &=& e_n + \int_{t_n}^{t_{n+1}} A\brac{s_1/\varepsilon} ds_1 e_n + \int_{t_n}^{t_{n+1}} A\brac{s_1/\varepsilon} \int_{t_n}^{s_1}  A\brac{s_2/\varepsilon} \brac{U(s_2)-U_n} ds_2 ds_1  \\ 
        && + \int_{t_n}^{t_{n+1}} A\brac{s_1/\varepsilon} \int_{t_n}^{s_1} \brac{g(U(s_2))-g_n}  ds_2 ds_1  + \int_{t_n}^{t_{n+1}} \brac{ g\brac{U(t_n)+h(s_1)} - g_n - \nabla \left. g \right|_n \tilde{h}(s_1) }  ds_1 \\
        &=& e_n + I+II+III+IV. 
    \end{eqnarray*}
Let us first consider the term $II$:
\begin{eqnarray}
II&=&    \int_{t_n}^{t_{n+1}} A\brac{s_1/\varepsilon} \int_{t_n}^{s_1}  A\brac{s_2/\varepsilon} \brac{U(s_2)-U_n} ds_2 ds_1 = \int_{t_n}^{t_{n+1}} A\brac{s_1/\varepsilon} \int_{t_n}^{s_1}  A\brac{s_2/\varepsilon} \brac{U(s_2)-U(t_n) + e_n} ds_2 ds_1 \nonumber\\
\label{term2}
    &=&\!\!\!\! \int_{t_n}^{t_{n+1}} \!\!\!\!A\brac{s_1/\varepsilon} \int_{t_n}^{s_1} \!\! A\brac{s_2/\varepsilon} \int_{t_n}^{s_2} \!\!\Big(\!{A\brac{s_3/\varepsilon} U(s_3) + g(U(s_3)) }\!\Big)  ds_3 ds_2 ds_1 
    \!+\! \int_{t_n}^{t_{n+1}} \!\!\!\!A\brac{s_1/\varepsilon} \int_{t_n}^{s_1}  \!\!A\brac{s_2/\varepsilon}  ds_2 ds_1 \; e_n.
\end{eqnarray}
From the assumptions on $A$ and $g$, we deduce $\|II\|\leq C\Delta t^3 + C\Delta t^2 \; e_n$. 
Similarly for the term $III$, we have $\|III\|\leq C\Delta t^2 \; e_n$. 
Considering the term $IV$, we have  
\begin{eqnarray*}
    IV &=& \int_{t_n}^{t_{n+1}} \brac{ g\brac{U(t_n)+h(s_1)} - g_n - \nabla \left. g \right|_n\tilde{h}(s_1) } ds_1  
    \\ 
    &=& \int_{t_n}^{t_{n+1}} \brac{ g\brac{U(t_n)+h(s_1)} - g\brac{U(t_n)} - (\nabla  g){(U(t_n))}\tilde{h}(s_1)  } ds_1 \\
    &&+ \int_{t_n}^{t_{n+1}} \brac{ g\brac{U(t_n)}-g_n +  \brac{ (\nabla  g){(U(t_n))} - \nabla \left. g \right|_n }\tilde{h}(s_1)  } ds_1. 
\end{eqnarray*}
From Lipshitz property on $g$ and $\nabla g$ and from the estimate $\|\tilde{h}(s_1)\|\leq C\Delta t$, we deduce that the 
second term on the right hand side is bounded 
by $C\Delta t e_n$. 
Upon Taylor expanding $g\brac{U(t_n)+h(s_1)}$ and considering the integral remainder, the first term on the right hand side of above equation becomes 
\begin{equation*}
    \int_{t_n}^{t_{n+1}} (\nabla  g){(U(t_n))} (h(s_1)-\tilde{h}(s_1))   ds_1  + \int_{t_n}^{t_{n+1}} \int_0^1 \brac{1-\tau} 
   ({\sf Hess}g){\brac{U(t_n)+\tau h(s_1)}}(h(s_1),h(s_1))
    d\tau ds_1  
\end{equation*}
where ${\sf Hess}g$ denotes the Hessian matrix of $g$ with respect to the variable $U$. Further, we have
\begin{eqnarray*}
\|(h-\tilde{h})(t)\| &=& \left\| \int_{t_n}^t A(s/\varepsilon)(U(s)-U_n) ds  + \int_{t_n}^t \Big[ g(U(s)) -g_n\Big] ds \right\|, 
\end{eqnarray*}
but $\|g(U(s)) -g_n\| \leq K\|U(s) - U_n\|\leq K (e_n + C\Delta t)$ hence $|(h-\tilde{h})(t)|\leq C\Delta t^2$. 
Further we also observe the following for $t_n<s_1<t_{n+1}$, with the notation $K_1=\norm{{\sf Hess}g({U})}$ 
\begin{eqnarray*}
    \norm{
    {\sf Hess}g({U})(h(s_1),h(s_1)) 
    } \leq K_1 \norm{h(s_1)^2} \leq K_1 \Delta t^2,   
\end{eqnarray*}  
so that we deduce $\|IV\|\leq C\Delta t^3 + C\Delta t \; e_n$. 
Then the norm of $e_{n+1}$ becomes:
$\norm{e_{n+1}} \leq (1+C\Delta t + C\Delta t^2) \norm{e_{n}} + C\Delta t^3$ 
and we conclude from the Gronwall lemma that 
$\norm{e_{n}} \leq C\Delta t^2$. 
\end{proof}

\begin{remark}
    Up to some technical and lengthy developments, it is possible to generalize the second order scheme presented above to arbitrary order of accuracy. 
\end{remark}

\subsection{Midpoint schemes for the nonlinear charged particles model}
\label{Section:Energy preserving scheme for avg model}
To handle the nonlinear models through midpoint schemes, fixed point strategies are usually required. To overcome 
this drawback, we introduce a scalar auxiliary variable following SAV techniques introduced in~\cite{sav2}.
More precisely, considering the  system \eqref{ODE(x,q)} and the  averaged model \eqref{sec4:averaged} 
which enjoys a Hamiltonian structure with the Hamiltonian \eqref{ham_norm}, 
our goal is to design uniformly accurate scheme for \eqref{ODE(x,q)}  which degenerates into an energy preserving scheme for the averaged model \eqref{sec4:averaged}.  Let us remark that contrary to the linear case, we do not have preservation of both Hamiltonians $H_1$ and $H_2$ (defined in \eqref{def_subh}).  
\subsubsection{SAV formulation of the averaged model}
In the nonlinear case, the averaged model can be written as 
\begin{equation}
\label{sec4:averaged}
    \begin{gathered}
    \dot{\bar{x}} = \bar{q} + \langle\theta\rangle J\bar{x}, \\
    \dot{\bar{q}} = -\nabla\phi(\bar{x}) + \langle\theta\rangle J\bar{q} +  \langle\theta^2\rangle J^2 \bar{x}, 
\end{gathered}
\end{equation}
and has a Hamiltonian structure where the Hamiltonian is given by \eqref{ham_norm}. 
Moreover, it is known that, except in some specific case or when using nonlinear solvers (e.g. averaged gradient method), 
there is no general numerical scheme preserving the Hamiltonian. 
We instead consider an extended ODE system based on a modification of the Hamiltonian 
through the SAV technique \cite{sav2}. Here, using the additional variable $\bar{r}(t) = \exp(\phi(\bar{x}(t))$, it is possible to reformulate the system~\eqref{sec4:averaged} into 
\begin{equation}
\label{sec4:averaged_sav}
    \begin{gathered}
    \dot{\bar{x}} = \bar{q} +   \langle\theta\rangle J\bar{x}, \\
    \dot{\bar{q}} = -\bar{b} + \langle\theta\rangle J\bar{q} +  \langle\theta^2\rangle J^2 \bar{x}, \\
    \dot{\bar{r}} = \bar{r} \bar{b}\cdot   \dot{\bar{x}},  
\end{gathered}
\end{equation}
where we introduced $\bar{b}(t)=\nabla \phi(\bar{x}(t))$. Other choices for the additional variable $\bar{r}$ are possible, the most used being $\bar{r}(t) = \sqrt{\phi(\bar{x}(t)) + C}$ with $C$ large enough.
This extended system enjoys an energy preservation property as shown in the next proposition. 
\begin{proposition}
\label{prop_sav_nl}
For the extended system~\eqref{sec4:averaged_sav}, the following quantity is preserved with time 
\begin{equation}
\label{Htilde}
\bar{H}(\bar{x},\bar{q},\bar{r}) = \frac{1}{2}|\bar{q}|^2 + \langle\theta\rangle \bar{q}\cdot J\bar{x} + \frac{1}{2}\langle\theta^2\rangle |\bar{x}|^2 + \log(\bar{r}). 
\end{equation}
\end{proposition}
\begin{proof}
Multiplying the first equation of~\eqref{sec4:averaged_sav} by $\langle\theta^2\rangle \bar{x}$, the second by $\bar{q}$ and combining them with the third equation gives 
\begin{eqnarray*}
\frac{d}{dt}\Big(\frac{1}{2}|\bar{q}|^2 +  \frac{1}{2}\langle\theta^2\rangle |\bar{x}|^2 + \log(\bar{r})\Big)&=&  \langle\theta^2\rangle \bar{x}\cdot \bar{q} - \bar{b}\cdot \bar{q} + \langle\theta^2\rangle J^2 \bar{x}\cdot \bar{q} +\bar{b}\cdot   \dot{ \bar{x}} \nonumber\\
&=&  - \bar{b}\cdot \Big( \dot{\bar{x}}  - \langle\theta\rangle J \bar{x}\Big) + \bar{b}\cdot \dot{\bar{x}} = \bar{b}\cdot  \langle\theta\rangle J\bar{x}. 
\end{eqnarray*}
The identity $-\bar{b}\cdot J\bar{x}=\frac{d}{dt}(\bar{q}\cdot J\bar{x})$ follows from 
\begin{eqnarray*}
\frac{d}{dt}(\bar{q}\cdot J\bar{x}) = \dot{\bar{q}}\cdot J\bar{x} + \bar{q}\cdot J\dot{\bar{x}} &=& -\bar{b}\cdot J\bar{x}+ \langle\theta\rangle J{\bar{q}}\cdot J\bar{x} +  \langle\theta\rangle \bar{q}\cdot J^2 \bar{x} =  -\bar{b}\cdot J\bar{x}\nonumber
\end{eqnarray*}
\end{proof}

\subsubsection{Numerical scheme for the SAV formulation of the averaged model}
\label{sec422:midpoint_nonlinear}
Next, we propose the following numerical scheme for solving the nonlinear averaged model~\eqref{sec4:averaged_sav}: 
\begin{equation}
\label{sav_nl_avg}
    \begin{gathered}
    (\bar{x}_{n+1} - \bar{x}_n) = \Delta t\bar{q}_{n+1/2} +  \Delta t\langle\theta\rangle J\bar{x}_{n+1/2}, \\
  (\bar{q}_{n+1} - \bar{q}_n)  =  - \Delta t \bar{b} + \Delta t\langle\theta\rangle J \bar{q}_{n+1/2} -  \Delta t\langle\theta^2\rangle  \bar{x}_{n+1/2}, \\
  \log(\bar{r}_{n+1}) - \log(\bar{r}_{n}) = \bar{b}\cdot  (\bar{x}_{n+1} - \bar{x}_n), 
\end{gathered}
\end{equation}
where $x_{n+1/2}=(x_{n+1}+x_n)/2$ and $q_{n+1/2}=(q_{n+1}+q_n)/2$. The approximation for the quantity $\bar b$ (which can be viewed as an approximation of  $(\int_{t_n}^{t_{n+1}}\nabla\phi(\bar{x}(s))ds)/\Delta t$) will be discussed later on and has no consequence on the following result: 
\begin{proposition}
\label{prop_sav_lin_energy}
The numerical solution to the scheme~\eqref{sav_nl_avg} preserves the Hamiltonian $\bar{H}$ given by~\eqref{Htilde}: $\bar{H}(\bar{x}_{n+1}, \bar{q}_{n+1}, \bar{r}_{n+1}) = \bar{H}(\bar{x}_{n}, \bar{q}_{n}, \bar{r}_n)$, $n\geq 0$. 
\end{proposition}
\begin{proof}
The proof follows the lines of the proof of Proposition \ref{prop_sav_nl}. Multiplying  the first equation 
of~\eqref{sav_nl_avg} by $\langle \theta^2 \rangle \bar{x}_{n+1/2}$, the second equation by  $\bar{q}_{n+1/2}$ and combining them 
with the third equation gives 
\begin{eqnarray}
\frac{1}{2} \langle\theta^2\rangle (|\bar{x}_{n+1}|^2 - |\bar{x}_{n}|^2) + \frac{1}{2}(|\bar{q}_{n+1}|^2 - |\bar{q}_{n}|^2) + \log(\bar{r}_{n+1}) - \log(\bar{r}_{n})  &=&\nonumber\\ 
&& \hspace{-10cm} \Delta t \langle\theta^2\rangle \bar{x}_{n+1/2} \cdot \bar{q}_{n+1/2} -\Delta t \bar{b}\cdot \bar{q}_{n+1/2} -  \Delta t \langle\theta^2\rangle \bar{x}_{n+1/2} \cdot \bar{q}_{n+1/2} +\bar{b}\cdot  (\bar{x}_{n+1} - \bar{x}_n)\nonumber\\
\label{sav_proof}
&& \hspace{-10.5cm} =-\bar{b}\cdot \Big((\bar{x}_{n+1} - \bar{x}_n) -  \Delta t \langle\theta\rangle J\bar{x}_{n+1/2}\Big) +\bar{b}\cdot  (\bar{x}_{n+1} - \bar{x}_n)=  \Delta t \langle\theta\rangle\bar{b}\cdot  J\bar{x}_{n+1/2}.  
\end{eqnarray}
It remains to check $-\Delta t \bar{b}\cdot  J\bar{x}_{n+1/2} = \bar{q}_{n+1}\cdot J\bar{x}_{n+1} - \bar{q}_{n}\cdot J\bar{x}_{n}$. \\
Multiplying the second equation of~\eqref{sav_nl_avg} by $J\bar{x}_{n+1/2}$ leads to 
$$
 (\bar{q}_{n+1} - \bar{q}_n) \cdot J\bar{x}_{n+1/2} = -\Delta t \bar{b}\cdot  J\bar{x}_{n+1/2} +  \Delta t \langle\theta\rangle J \bar{q}_{n+1/2} \cdot J\bar{x}_{n+1/2} =-\Delta t \bar{b}\cdot  J\bar{x}_{n+1/2} +  \Delta t \langle\theta\rangle  \bar{q}_{n+1/2} \cdot \bar{x}_{n+1/2}. 
$$
Considering $\bar{q}_{n+1/2}\cdot J$ multiplied by the first equation of~\eqref{sav_nl_avg} leads to 
\begin{eqnarray*} 
\bar{q}_{n+1/2} \cdot  J(\bar{x}_{n+1} - \bar{x}_n) &=& \Delta t \langle\theta\rangle \bar{q}_{n+1/2} \cdot J^2 \bar{x}_{n+1/2} = -\Delta t \langle\theta\rangle \bar{q}_{n+1/2} \cdot \bar{x}_{n+1/2}. 
\end{eqnarray*}
Adding the two last equalities gives 
$$
(\bar{q}_{n+1} - \bar{q}_n) \cdot J\bar{x}_{n+1/2} + \bar{q}_{n+1/2} \cdot  J(\bar{x}_{n+1} - \bar{x}_n) = -\Delta t \bar{b}\cdot  J\bar{x}_{n+1/2},  
$$
and using~\eqref{relation_discrete}, one gets $
-\Delta t \bar{b}\cdot  J\bar{x}_{n+1/2} = \bar{q}_{n+1} \cdot J\bar{x}_{n+1} - \bar{q}_n\cdot J\bar{x}_n$ which enables to conclude with \eqref{sav_proof}. 
\end{proof}

We end this part by discussing the approximation of $\bar{b}$. The usual choice is 
$\bar{b} \approx \nabla \phi(\bar{x}(t_{n+1/2}))$ so that second order accuracy is reached through extrapolation techniques \cite{sav2}: 
\begin{equation}
\label{bbar_extrap}
\bar{b} \approx \nabla \phi(\bar{x}(t_{n+1/2})) \approx -\frac{1}{2} \nabla \phi(\bar{x}_{n-1}) +\frac{3}{2}\nabla \phi(\bar{x}_{n}). 
\end{equation}
Even if this choice enables to get a linearly implicit scheme, it requires the knowledge of $x_{n-1}$. 
Here, another technique is proposed based on the above schemes derived from integration in time. Indeed, 
we start from 
\begin{equation}
\label{b_bar}
\bar{b}= \frac{1}{\Delta t}\int_{t_n}^{t_{n+1}} \nabla \phi(\bar{x}(s))ds. 
\end{equation}
But from the integration on $t\in [{t_n}, s]$ 
of the first equation of \eqref{sec4:averaged_sav}, we get 
$$
\bar{x}(s)= \bar{x}(t_n) + \bar{h}(s)  \;\; \mbox{ with }  \;\; \bar{h}(s) = \int_{t_n}^s \bar{q}(s_1) ds_1 +  \int_{t_n}^s \langle \theta \rangle J \bar{x}(s_1) ds_1.  
$$
We can now approximate $\bar{h}$ by $\tilde{\bar{h}}=(s-{t_n}) \bar{q}(t_n) +  (s-{t_n}) \langle \theta \rangle J \bar{x}(t_n)$ 
so that a choice for $\bar{b}$ is obtained by inserting this last approximation in a Taylor expansion for $\nabla \phi$ in \eqref{b_bar}
\begin{eqnarray}
\bar{b} &=& \frac{1}{\Delta t}\int_{t_n}^{t_{n+1}} \nabla \phi(\bar{x}(s))ds = \frac{1}{\Delta t}\int_{t_n}^{t_{n+1}} \nabla \phi(\bar{x}(t_n)+\bar{h}(s))ds \nonumber\\ 
\label{bar_b}
&\approx & \nabla\phi(\bar{x}(t_n)) \!+ \!\nabla^2\phi(\bar{x}(t_n))  \frac{1}{\Delta t}\int_{t_n}^{t_{n+1}} \tilde{\bar{h}}(s) ds \nonumber\\
&\approx &  \nabla\phi(\bar{x}_n) + \nabla^2\phi (\bar{x}_n)\frac{\Delta t}{2}\Big(\bar{q}_n + \langle \theta \rangle J \bar{x}_n\Big),  
\end{eqnarray}
which can be proved to be a second order (third locally) approximation of $(\int_{t_n}^{t_{n+1}} \nabla \phi(\bar{x}(s))ds)/\Delta t$ and only requires the knowledge of $\bar{x}_n$ to update the numerical unknown  $(\bar{x}_{n+1},\bar{q}_{n+1},\bar{r}_{n+1})$.

\subsubsection{Uniformly accurate and structure preserving SAV-schemes  for the charged particle model}
\label{sec423:midpoint_nonlinear}
The goal of this part is to propose a numerical scheme for the charged particle model~\eqref{ODE(x,q)} which is second order uniformly accurate and which degenerates when $\varepsilon\to 0$ towards 
the SAV-midpoint scheme~\eqref{sav_nl_avg} (which preserves the 
energy thanks to Proposition \ref{prop_sav_lin_energy}). First, 
we derive the SAV formulation for~\eqref{ODE(x,q)}, inspired by the SAV formulation of the averaged model
\begin{equation}
\label{sec4:sav}
    \begin{gathered}
    \dot{{x}} = {q} +   \theta J{x}, \\
    \dot{{q}} = -b + \theta J{q} + \theta^2 J^2 {x}, \\
    \dot{{r}} =  r b\cdot   \dot{{x}}.  
\end{gathered}
\end{equation}
To derive a second order UA scheme, let us integrate 
the first two equations of \eqref{sec4:sav} on $t\in[t_n, t_{n+1}]$ to get 
\begin{eqnarray}
\label{x_ex}
x(t_{n+1})-x(t_n) &=& \int_{t_n}^{t_{n+1}} q(t)dt + \int_{t_n}^{t_{n+1}} \theta(t/\varepsilon) Jx(t)dt, \\
\label{q_ex}
q(t_{n+1})-q(t_n) &=& -\int_{t_n}^{t_{n+1}} b(t)dt + \int_{t_n}^{t_{n+1}} \theta(t/\varepsilon) Jq(t)dt+ \int_{t_n}^{t_{n+1}} \theta^2(t/\varepsilon) J^2 x(t)dt,
\end{eqnarray}
but following the calculations performed in \eqref{exact_midpoint}, we get 
\begin{eqnarray*}
x(t)&=&\frac{1}{2}(x(t_{n+1}+x(t_n))+ \frac{1}{2}\Big[\int_{t_n}^{t}(q(s)+\theta(s/\varepsilon)Jx(s))ds-\int_{t}^{t_{n+1}}(q(s)+\theta(s/\varepsilon)Jx(s))ds\Big] \nonumber\\
q(t)&=&\frac{1}{2}(q(t_{n+1}+q(t_n))-  \frac{1}{2} \Big[\int_{t_n}^t b(s)ds - \int_t^{t_{n+1}}b(s)ds \Big] \nonumber\\
&&+\frac{1}{2}\Big[\int_{t_n}^{t} (\theta(s/\varepsilon) Jq(s)+\theta^2(s/\varepsilon)J^2x(s))ds- \int_t^{t_{n+1}}(\theta(s/\varepsilon) Jq(s)+\theta^2(s/\varepsilon)J^2x(s))ds\Big].  
\end{eqnarray*}
Inserting these expressions into \eqref{x_ex} and \eqref{q_ex}, and using midpoint approximation of $(x(s), q(s))\approx (x_{n+1/2}, q_{n+1/2})$ leads to  
\begin{eqnarray}
\label{x_app0}
x_{n+1}-x_n &=& \!\!(\mathcal{B}_{1,1} + \frac{1}{2}\mathcal{A}_{1,1}) x_{n+1/2} + (\mathcal{B}_{1,2} + \frac{1}{2}\mathcal{A}_{1,2})  {q}_{n+1/2} -\frac{1}{2}\int_{t_n}^{t_{n+1}}\!\!\Big[  \int_{t_n}^t \!\!\! b(s)ds -\! \int_{t}^{t_{n+1}}\!\!\! b(s) ds \Big]dt, \\
q_{n+1}-q_n &=& -\int_{t_n}^{t_{n+1}}b(s)ds + (\mathcal{B}_{2,1} +\frac{1}{2} \mathcal{A}_{2,1})  x_{n+1/2} +(\mathcal{B}_{2,2} + \frac{1}{2} \mathcal{A}_{2,2})   q_{n+1/2} \nonumber\\
\label{q_app0}
&&-\frac{1}{2}\int_{t_n}^{t_{n+1}}\theta(t/\varepsilon)J\Big[  \int_{t_n}^t b(s)ds - \int_t^{t_{n+1}} b(s)ds\Big]dt, 
\end{eqnarray}
where the $2$x$2$ block matrices are given by~\eqref{a11},~\eqref{a12},~\eqref{a21},~\eqref{a22} 
and the $4$x$4$ matrix $\mathcal{B}$  is given by \eqref{midpoint_B}. The $b$ terms with double integrals can be approximated at $s=t_n$ 
(the last term in \eqref{x_app0} vanishes)  
whereas the term $-\int_{t_n}^{t_{n+1}}b(s)ds$ will be discussed after. 
We thus consider the following numerical scheme for~\eqref{sav_nl_avg} 
\begin{equation}
    \begin{gathered}
    \label{sav_nl_mid}
    x_{n+1} - x_n = (\mathcal{B}_{1,1} + \frac{1}{2}\mathcal{A}_{1,1}) x_{n+1/2} + (\mathcal{B}_{1,2} + \frac{1}{2}\mathcal{A}_{1,2})  {q}_{n+1/2}, \\
  q_{n+1} - q_n  =  -\int_{t_n}^{t_{n+1}} {b}(s)ds + (\mathcal{B}_{2,1} +\frac{1}{2} \mathcal{A}_{2,1})  x_{n+1/2} +(\mathcal{B}_{2,2} + \frac{1}{2} \mathcal{A}_{2,2})   q_{n+1/2}-\int_{t_n}^{t_{n+1}}\theta(t/\varepsilon)(t-t_{n+1/2})dt J b_n, \\
  \log(r_{n+1}) - \log(r_{n}) = \int_{t_n}^{t_{n+1}} {b}(s)ds \cdot  (x_{n+1} - x_n), 
\end{gathered}
\end{equation} 
where $b_n=\nabla \phi(x_n)$ and $t_{n+1/2}=(t_n+t_{n+1})/2$. 
Let us remark that the only requirement for the approximation of $\int_{t_n}^{t_{n+1}} {b}(s)ds$ is to be second order accurate approximation of 
$\int_{t_n}^{t_{n+1}} \nabla\phi(x(s))ds$. Below, we present two ways to approximate $\int_{t_n}^{t_{n+1}} {b}(s)ds$. 
\paragraph{Choice 1} 
The first choice proposed by~\cite{sav2} is 
\begin{equation}
\label{bad_b}
\frac{1}{\Delta t}\int_{t_n}^{t_{n+1}} {b}(s)ds \approx \nabla \phi(x_{n+1/2}) = -\frac{1 }{2}\nabla\phi(x_{n-1})+\frac{3}{2}\nabla\phi(x_{n}).
\end{equation}

\paragraph{Choice 2} 
The second choice follows the strategy we used to construct uniformly accurate scheme in the nonlinear case. 
First, we have 
$$
 \frac{1}{\Delta t}\int_{t_n}^{t_{n+1}} \nabla \phi(x(s)) ds = \frac{1}{\Delta t}\int_{t_n}^{t_{n+1}}\nabla \phi(x(t_n) + h(s)) ds,  
$$
where $h(s)= \int_{t_n}^s q(s_1) ds_1 +\int_{t_n}^s \theta(s_1/\varepsilon) J x(s_1) ds_1$ is obtained by integrating the equation on $x$ in \eqref{sec4:sav}.
In view of a second order numerical scheme, 
we approximate $h(s)$ by $\tilde{h}(s)$ given by  
\begin{equation}
\label{htilde}
\tilde{h}(s) = (s-t_n) q(t_n) + \int_{t_n}^s \theta(s_1/\varepsilon) J  ds_1 x(t_n). 
\end{equation}
Since $\phi$ is assumed to be smooth, we perform a Taylor expansion 
\begin{eqnarray*}
\int_{t_n}^{t_{n+1}}\nabla \phi(x(t_n) + {h}(s)) ds &\approx& \int_{t_n}^{t_{n+1}}\nabla\phi(x(t_n) + \tilde{h}(s)) ds\approx \Delta t \nabla \phi(x(t_n))  +\nabla^2 \phi(x(t_n))\int_{t_n}^{t_{n+1}} \tilde{h}(s) ds \nonumber\\
&\approx& \Delta t\nabla \phi(x(t_n))  +\nabla^2 \phi(x(t_n))\int_{t_n}^{t_{n+1}} \Big[(s-t_n) q(t_n) + \int_{t_n}^s \theta(s_1/\varepsilon) J  ds_1 x(t_n) \Big]ds \nonumber\\
&\approx& \Delta t\nabla \phi(x(t_n))  +\nabla^2 \phi(x(t_n))\Big[\frac{\Delta t^2}{2} q(t_n) + \int_{t_n}^{t_{n+1}} \int_{t_n}^s \theta(s_1/\varepsilon) J  ds_1 ds x(t_n) \Big]. 
\end{eqnarray*}
As a conclusion, $\frac{1}{\Delta t}\int_{t_n}^{t_{n+1}} {b}(s)ds$ is chosen as follows 
\begin{equation}
\label{good_b}
\frac{1}{\Delta t}\int_{t_n}^{t_{n+1}} {b}(s)ds = \nabla \phi(x_n) +\frac{\Delta t}{2} \nabla^2 \phi(x_n) q_n + \frac{1}{\Delta t}\nabla^2 \phi(x_n) \int_{t_n}^{t_{n+1}} \int_{t_n}^s \theta(s_1/\varepsilon) J  ds_1 ds \, x_n. 
\end{equation}

We end this part by identifying the limit of \eqref{sav_nl_mid} when $\varepsilon\to 0$ with the energy preserving scheme \eqref{sav_nl_avg}. 
\begin{proposition}
The numerical scheme for charged particle model~\eqref{sav_nl_mid} with the approximation  given by~\eqref{good_b} degenerates when $\varepsilon \to 0$ 
to the energy preserving numerical scheme for averaged model~\eqref{sav_nl_avg} with $\bar{b}$ given by~\eqref{bar_b}. 
\end{proposition}

\begin{proof}
The only point to check is the limit of \eqref{good_b} 
when $\varepsilon\to 0$. From the definition of  given by~\eqref{good_b}, we get (assuming $x_n\to \bar{x}_n$ and $q_n\to \bar{q}_n$ when $\varepsilon\to 0$) 
\begin{eqnarray*}
\frac{1}{\Delta t}\int_{t_n}^{t_{n+1}} {b}(s)ds &=&  \nabla \phi(\bar{x}_n) + \nabla^2 \phi(\bar{x}_n) \Big[ \frac{\Delta t}{2} \bar{q}_n + \frac{1}{\Delta t}\int_{t_n}^{t_{n+1}}(s-t_n)ds  \langle \theta \rangle J \bar{x}_n\Big] + {\cal O}(\varepsilon) \nonumber\\
&=&  \nabla \phi(\bar{x}_n) + \nabla^2 \phi(\bar{x}_n)  \frac{\Delta t}{2} \Big[ \bar{q}_n +   \langle \theta \rangle J \bar{x}_n\Big] + {\cal O}(\varepsilon)  =\bar{b}+ {\cal O}(\varepsilon) , 
\end{eqnarray*}
which is the definition \eqref{bar_b} of $\bar{b}$ given for the SAV averaged model.  
\end{proof}

\section{Particle-In-Cell}
\label{pic}

In this section, we are going to extend the numerical methods studied in previous sections for solving the Vlasov-Poisson equation (\ref{vlasov}) with a fast oscillating magnetic field. Indeed, using \cite{Bostan2012} and based on the two-scale convergence \cite{allaire}, the solution to \eqref{vlasov} can be proved to converge 
to a Vlasov equation whose characteristics are given by \eqref{sec4:averaged_sav}. Hence, we extend our schemes developed in the ODE framework to the Vlasov case by  considering the Particle-in-Cell (PIC) discretization. The PIC discretization approximates the unknown distribution $f^\varepsilon(t,x,v)$ of (\ref{vlasov}) by a sum of Dirac masses located at $(x_k(t), v_k(t))$ with weight $\omega_k>0$ for $k=1,\ldots,N_p$ and $N_p\in\mathbb{N}$ as
\begin{equation}
\label{dirac}
f^\varepsilon_{p}(t,x,v)=\sum_{k=1}^{N_p}\omega_k\delta(x-x_k(t))\delta(v-v_k(t)), \quad t\geq0,\ x, v\in\mathbb{R}^2.
\end{equation}
The weights $\omega_k$ and initial values of the particles $x_{k,0},v_{k,0}$ for $k=0,\ldots,N_p$ are prescribed according to the  initial condition $f^\varepsilon(t=0, x, v)=f^{\text{in}}(x,v)$ in (\ref{vlasov}). Considering uniform weights, we have 
$$
\omega_k=\frac{1}{N_p}\int_{\mathbb{R}^2\times\mathbb{R}^2}f^{\text{in}}(x,v)dx dv,\quad k=1,\ldots,N_p.
$$
Standard sampling techniques like the Monte Carlo type rejection sampling method found in standard textbooks~\cite{SonnendruckerBook} can be employed to obtain the initial data for position and velocity of the particles  $x_{k,0},v_{k,0}$ for $k=1,\ldots,N_p$, 
. 
However, when $f^{\text{in}}$ takes the form $f^{\text{in}}(x,v)=\chi(x)M(v)$ with $\chi$ periodic and $M$ Maxwellian (which is a widely used case in plasma physics), specific deterministic techniques can be employed using the inversion of cumulative distribution function (see~\cite{Birdsall}). 

The dynamics of macro-particles is then given by 
the characteristic equations  (\ref{charact}) which are coupled to the Poisson equation 
$$
-\Delta_x \phi^\varepsilon(t, x) = \nabla_x \cdot E^\varepsilon(t, x) = \rho^\varepsilon(t, x) - 1, \;\; \mbox{ with } \rho^\varepsilon(t, x) = \int_{\mathbb{R}^2} f^\varepsilon(t, x, v) dv, 
$$
through the potential of electric field as $E^\varepsilon = - \nabla_x \phi^\varepsilon$.
From the positions $\{x_k(t)\}_{k=1,\ldots,N_p}$ of the particles at time $t>0$, the density $\rho^\varepsilon$ is approximated by $\rho_p$ given by  
$$
\rho^\varepsilon(t,x) \approx \rho_p^\varepsilon(t,x)=\sum_{k=1}^{N_p}\omega_\ell \delta(x-x_k(t))\quad  x\in\mathbb{R}^2. 
$$
This enables to evaluate the density on a uniform grid so that the Poisson equation can be also solved on a mesh grid of $x$ in $\mathbb{R}^2$ to get $E_p^\varepsilon(t,x)\approx E^\varepsilon(t,x)$ 
Finally, an interpolation of $E_p^\varepsilon(t,x)$ is required at each particle position in the pusher. 
In practice, the Dirac function $\delta(x)$ is approximated by  B-spline function $S^m(x)$ ($m\in\mathbb{N}$) given by 
\cite{SonnendruckerBook}:
$$
S^0(x):=\left\{\begin{cases}&\frac{1}{\Delta x},\qquad |x|\leq\frac{\Delta x}{2},\\
&0,\qquad\quad \mbox{else},\end{cases}\right.\qquad S^m(x):=\frac{1}{\Delta x}\int_{x-\frac{\Delta x}{2}}^{x+\frac{\Delta x}{2}}S^{m-1}(y)dy,\quad m\geq 1.
$$
The case in two dimensions is obtained by tensor product. 

Under the PIC discretization (\ref{dirac}), the characteristic equations  of \eqref{vlasov}
 read for $k=1,\ldots,N_p$
\begin{subequations}
\label{charact}
\begin{align}
   & \dot{x}_k(t)=v_k(t), \label{characta}\\
&\dot{v}_k(t)=E(t,x_k(t))+\frac{B}{2\varepsilon}\theta'(t/\varepsilon) J x_k(t) + B \theta(t/\varepsilon) Jv_k(t),\quad t>0, \\
   & x_k(0)=x_{k,0},\quad v_k(0)=v_{k,0}, 
\end{align}
\end{subequations}
where $B$ denotes the magnetic field amplitude (the magnetic field ${\bf B}=\theta(t/\varepsilon)(0, 0, B)$ is considered).  
Under suitable condition on  $f^\text{in}$ (see~\cite{Bostan2012}), we will consider the Vlasov equation on $g^\varepsilon(t,x,q) = f^{\varepsilon}(t, x, v)$ 
using the change of variables $q=v-B/2 \; \theta(t/\varepsilon) Jx$. Hence,  we consider the following ODE system (for $k=1, \dots, N_p$)   
\begin{subequations}
\label{charact_xq}
\begin{align}
   & \dot{x}_k(t)=q_k(t) + \frac{B}{2} \theta(t/\varepsilon) Jx_k(t), \label{characta_xq}\\
  &\dot{q}_k(t)= E(t,x_k(t)) + \frac{B}{2}\theta(t/\varepsilon) J q_k(t) + \frac{B^2}{4}\theta^2(t/\varepsilon) J^2 x_k(t),\quad t>0, \label{charactb_xq}\\
   & x_k(0)=x_{k,0},\quad q_k(0)=v_{k,0}-\frac{B}{2} \; \theta(t/\varepsilon) Jx_{k,0}, 
\end{align}
\end{subequations}
Hence, the PIC method is applied to 
 $g^\varepsilon(t,x,q) \approx g_p(t, x, q) = \sum_{k=1}^{N_p} \omega_k \delta(x-x_k(t))\delta(q-q_k(t))$ where $x_k, q_k$ are solution to~\eqref{charact_xq} for which we developed 
uniformly accurate numerical schemes in the previous section.

\section{Numerical results}
\label{sec:numerical_results}
In this section, we present some numerical results to illustrate the  explicit and midpoint uniformly accurate methods for the  systems~\eqref{First order exact_ODE(x,q)} and~\eqref{averaged_ODE}. We also present the numerical results for Vlasov-Poisson equation by employing Particle-In-Cell method along with our uniformly accurate methods for ODE systems. 

\subsection{ODE systems} 
\begin{figure}
    \centering
    \subfloat[$L^2$ error of $x_1,x_2$ as a function of $\Delta t$]{\label{error_1_x_dt}\includegraphics[width=0.4\textwidth]{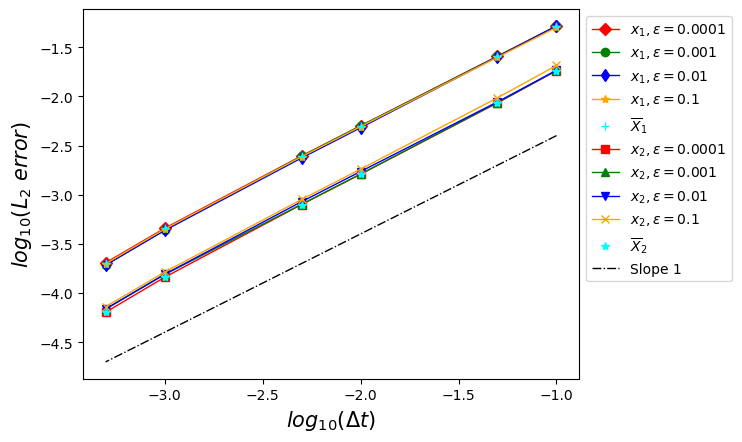}}
    \subfloat[$L^2$ error of $x_1,x_2$ as a function of $\varepsilon$]{\label{error_1_x_epsilon}\includegraphics[width=0.4\textwidth]{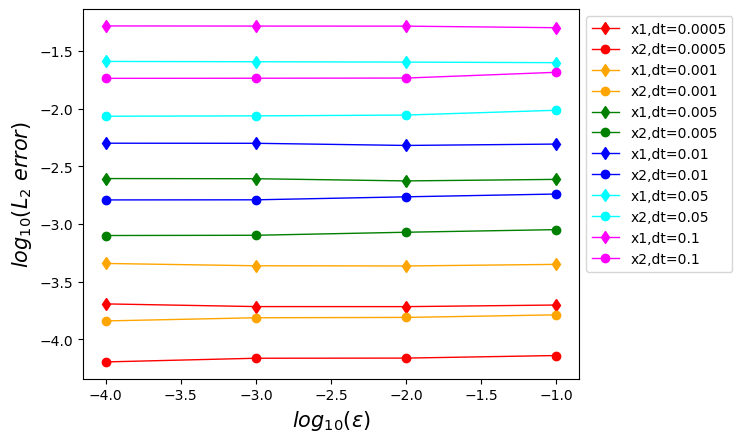}}
    
    \subfloat[$L^2$ error of $q_1,q_2$ as a function of $\Delta t$]{\label{error_1_q_dt}\includegraphics[width=0.4\textwidth]{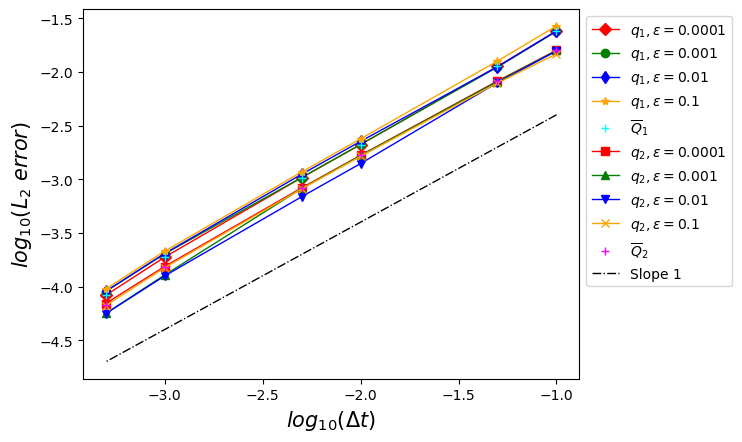}}
    \subfloat[$L^2$ error of $q_1,q_2$ as a function of $\varepsilon$]{\label{error_1_q_epsilon}\includegraphics[width=0.4\textwidth]{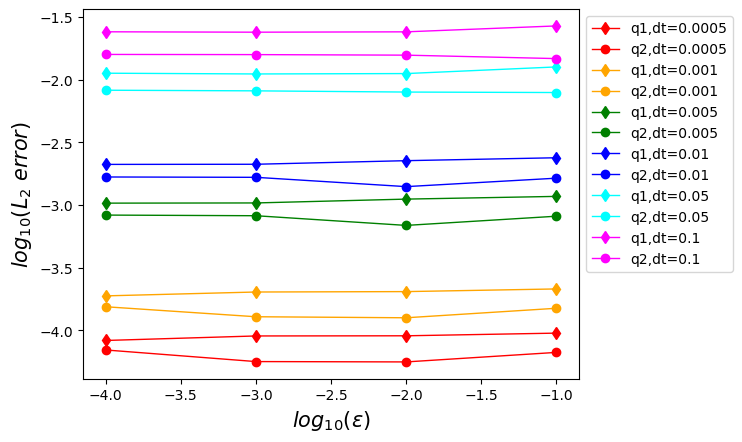}}
    \caption{\centering Linear case, first order explicit scheme: First order uniform accuracy with $\Delta t$.}
    \label{error_1}
\end{figure}

\begin{figure}
    \centering
    \subfloat[error with $\Delta t$ -  midpoint]{\label{error_ep_1_x_dt}\includegraphics[width=0.3\textwidth]{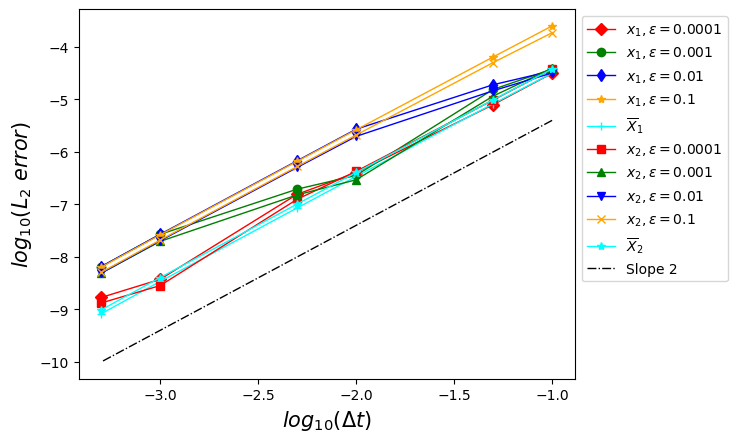}}
    \subfloat[error with $\Delta t$ - UA midpoint]{\label{error_ep_2_x_dt}\includegraphics[width=0.3\textwidth]{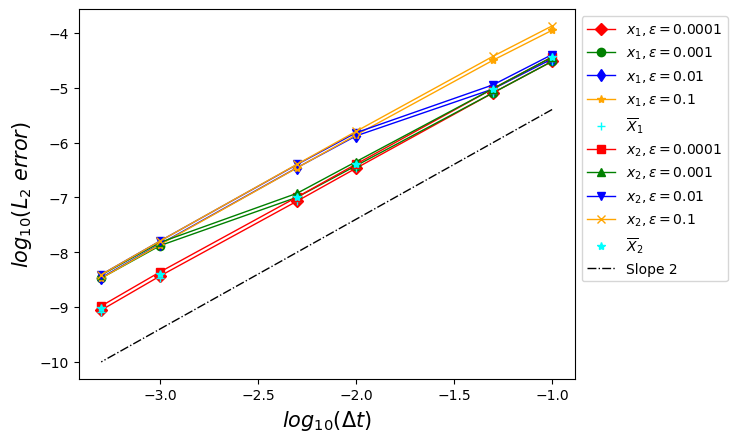}}
    \subfloat[error with $\varepsilon$ - UA midpoint]{\label{error_ep_2_x_epsilon}\includegraphics[width=0.3\textwidth]{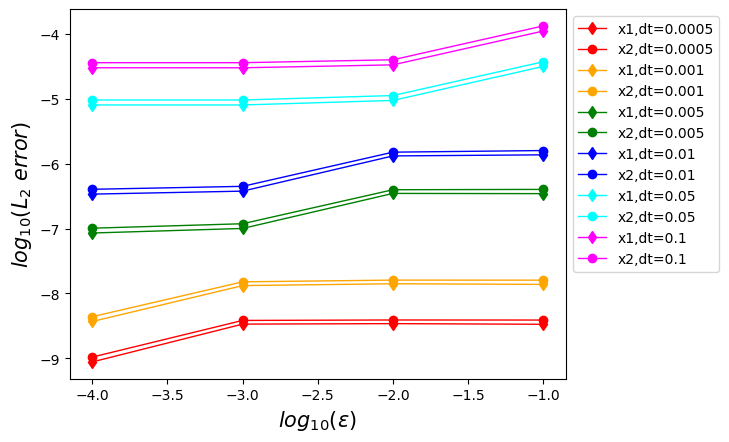}}

    \caption{\centering Linear case,  midpoint and UA midpoint schemes: Second order uniform accuracy with $\Delta t$.}
    \label{error_2_EP}
\end{figure}

\begin{figure}
    \centering
    \subfloat[$L^2$ error of $y$ as a function of $\Delta t$]{\label{error_order4_dt}\includegraphics[width=0.4\textwidth]{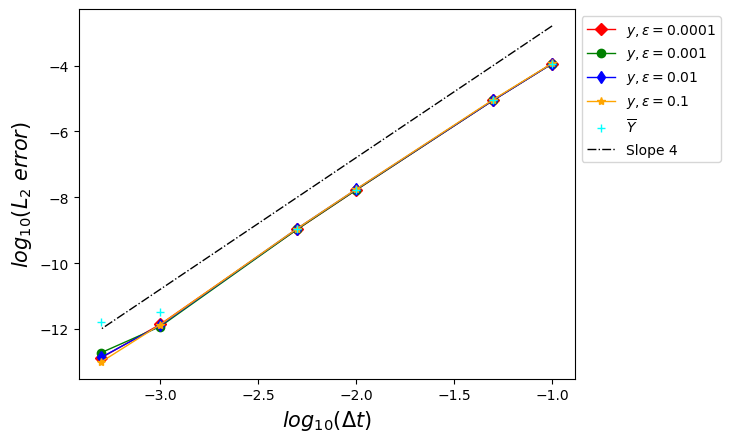}}
    \subfloat[$L^2$ error of $y$ as a function of  $\varepsilon$]{\label{error_order4_epsilon}\includegraphics[width=0.4\textwidth]{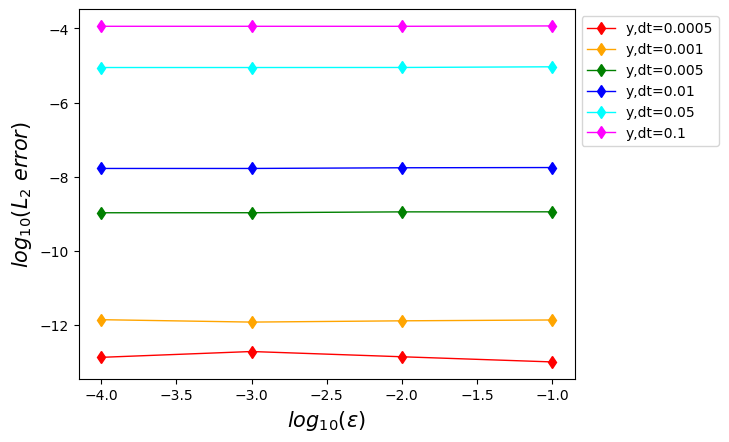}}
    \caption{\centering Linear scalar case, fourth order explicit scheme: uniform accuracy with $\Delta t$.}
    \label{error_4}
\end{figure}

The properties of numerical schemes for the ODE systems are validated in this subsection. 
\paragraph{UA property\\} 
The uniform accuracy (UA) properties of both explicit and energy preserving numerical schemes are presented in the linear case $\dot{U}(t)=A(t/\varepsilon)U(t), \; U(t=0)=U_0, t\in [0, T], T=1$. In \Cref{error_1,error_2_EP,error_4}, we plot the $L^2$ error 
to illustrate the order of (uniform) accuracy. More precisely, in 
Figure \ref{error_1}, the first order explicit UA scheme is considered with $U(t)=(x(t), q(t))\in\mathbb{R}^4$, $\theta(s)=1+\cos(s)$ in the matrix $A\in {\mathcal M}_{4,4}(\mathbb{R})$ given by 
\begin{equation}
\label{A_numerical_section}
A \brac{t/\varepsilon} = \begin{bmatrix} \theta\brac{t/\varepsilon} J & 1 \\ \theta^2\brac{t/\varepsilon}  J^2 & \theta\brac{t/\varepsilon} J \end{bmatrix}.  
\end{equation}
In  Figure \ref{error_2_EP}, the standard midpoint and UA midpoint schemes are considered with the same ODE as before except $\theta(s)=\cos(s)$. It is observed that second order accuracy is obtained for both schemes.
Finally, fourth order UA scheme is considered in Figure \ref{error_4} 
for the scalar model $\dot{U}(t)=[2+0.5\cos^2(t/\varepsilon)]U(t)$. 
For all the UA schemes studied in \Cref{error_1,error_2_EP,error_4}, 
the uniform accuracy is recovered. 

\begin{figure}
    \centering
    \subfloat[$L^2$ error of $x_1,x_2$ as a function of $\Delta t$]{\label{error_2_x_dt}\includegraphics[width=0.4\textwidth]{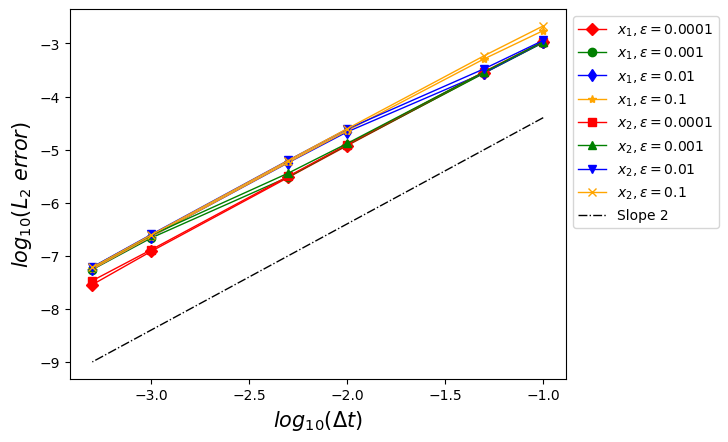}}
    \subfloat[$L^2$ error of $x_1,x_2$ as a function of $\varepsilon$]{\label{error_2_x_epsilonNL}\includegraphics[width=0.4\textwidth]{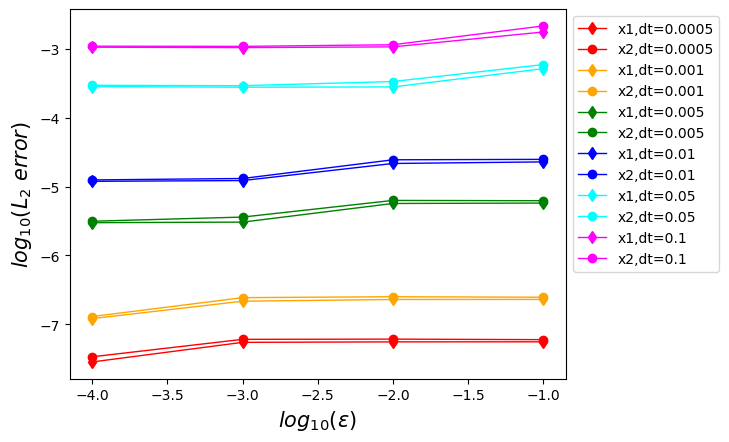}}

    \caption{\centering Nonlinear case, second order explicit scheme: Second order uniform accuracy with $\Delta t$.}
    \label{error_2_NL}
\end{figure}

\begin{figure}
    \centering
    \subfloat[$L^2$ error of $x_1,x_2$ as a function of $\Delta t$]{\label{error_2_x_dt_NLEP_badb_2}\includegraphics[width=0.4\textwidth]{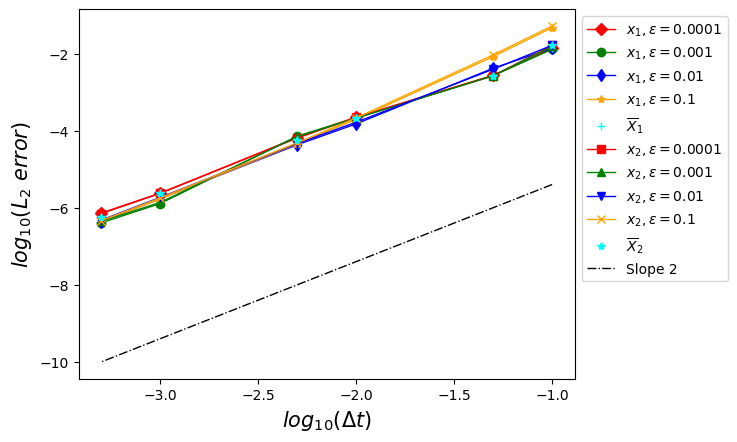}}
    \subfloat[$L^2$ error of $x_1,x_2$ as a function of $\varepsilon$]{\label{error_2_x_epsilon_NLEP_badb_2}\includegraphics[width=0.4\textwidth]{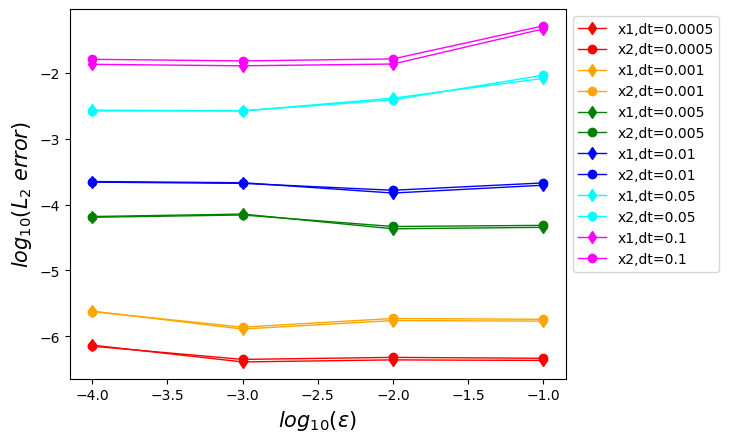}}
    \caption{\centering Nonlinear case, midpoint SAV-scheme (first choice for $b$): Second order uniform accuracy with $\Delta t$. }
    \label{error_2_NL_EP_badb}
\end{figure}

\begin{figure}
    \centering
    \subfloat[$L^2$ error of $x_1,x_2$ as a function of $\Delta t$]{\label{error_2_x_dt_NLEP}\includegraphics[width=0.4\textwidth]{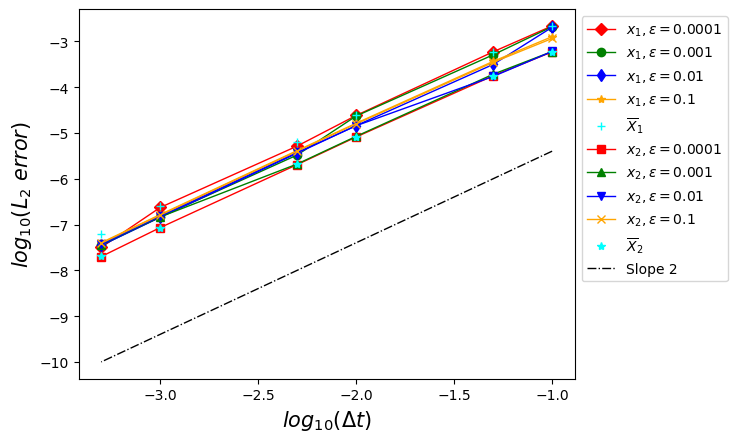}}
    \subfloat[$L^2$ error of $x_1,x_2$ as a function of $\varepsilon$]{\label{error_2_x_epsilon_NLEP}\includegraphics[width=0.4\textwidth]{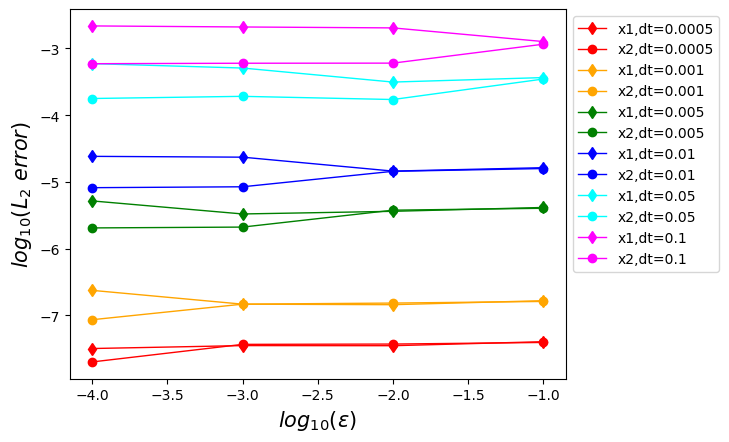}}
    \caption{\centering Nonlinear case, midpoint SAV-scheme $\!\!\!$ (second choice for $\!b$): Second order uniform accuracy with $\Delta t$.}
    \label{error_2_NL_EP}
\end{figure}

In Figures \ref{error_2_NL},\ref{error_2_NL_EP_badb},\ref{error_2_NL_EP}, 
the nonlinear case is considered: $\dot{U}(t)=A(t/\varepsilon)U(t)+g(U(t)), \; U(t=0)=U_0, t\in [0, T], T=1$, 
with $U(t)=(x(t), q(t))\in\mathbb{R}^4$ and $g(U(t))=(0, 0, \cos{x_1}\sin{x_2}+x_1+x_1^3, \sin{x_1}\cos{x_2}+x_2+x_2^3)^T\in\mathbb{R}^4$. Finally, the 
matrix $A$ given by~\eqref{A_numerical_section} is considered with $\theta(s)=\cos(s)$. 
Figure \ref{error_2_NL} depicts the results obtained by the second order explicit UA scheme \eqref{Second order_Enon0_ODE_(x,q)}. The expected second order uniform accuray is recovered.  Figures \ref{error_2_NL_EP_badb},\ref{error_2_NL_EP} correspond to midpoint second order scheme for the SAV reformulated system given by \eqref{sav_nl_mid}. We compare the two choices of $b$: in  
\Cref{error_2_NL_EP_badb}, the results are obtained with the first choice for the approximation of $\int_{t_n}^{t_{n+1}}{b}(s)ds$  given by~\eqref{bad_b} whereas  \Cref{error_2_NL_EP}  involves the approximation given by~\eqref{good_b}. It is observed that both choices reach the second order accuracy. 
Further, we also observe that the numerical scheme \eqref{sav_nl_avg} for the asymptotic model with both choices of $\bar{b}$ (given by \eqref{bbar_extrap} or \eqref{b_bar} preserve energy up to machine accuracy for long time (figures not shown).

\paragraph{Geometric property of the averaged model\\} 

In this part, we present the phase plots and frequency spectra obtained by energy preserving methods for the approximation of the following models: $\dot{U} = A(t/\varepsilon) U$,  $A$ is given by \eqref{A_numerical_section} with $\theta$ replaced by $B \theta$, $B\in\mathbb{R}$ being the magnitude of magnetic field. When $\varepsilon$ goes to zero, the asymptotic model is determined by the eigenvalues of $\langle A\rangle$ which can be computed for a given function $\theta$. 

\begin{figure}
    \centering
    \subfloat[$x_1(t)$ and $\bar{X}_1(t)$]{\label{x1t_alp1eps01}\includegraphics[width=0.3\textwidth]{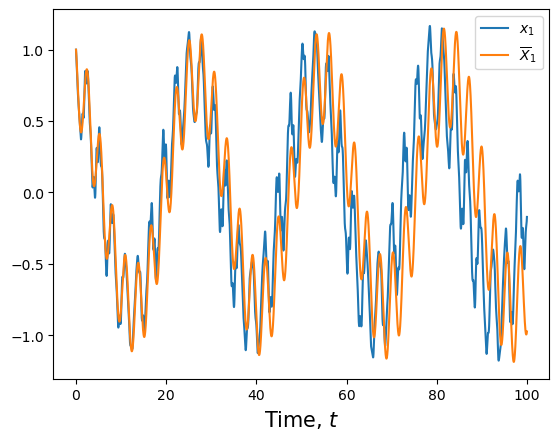}}
    \subfloat[Phase plot of $(x_1,x_2)$]{\label{x1x2_alp1eps01}\includegraphics[width=0.3\textwidth]{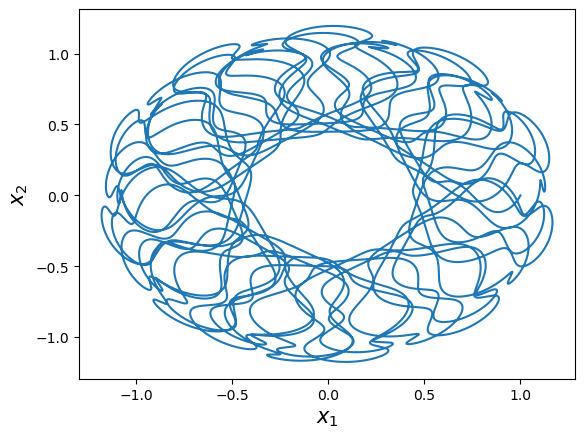}}
    \subfloat[Phase plot of $(\bar{X}_1,\bar{X}_2)$]{\label{x1bx2b_alp1eps01}\includegraphics[width=0.3\textwidth]{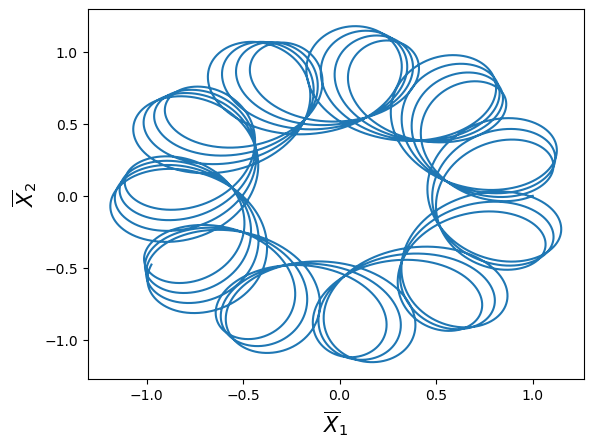}}

    \subfloat[Frequencies of $x_1$]{\label{FFTx1_alp1eps01}\includegraphics[width=0.4\textwidth]{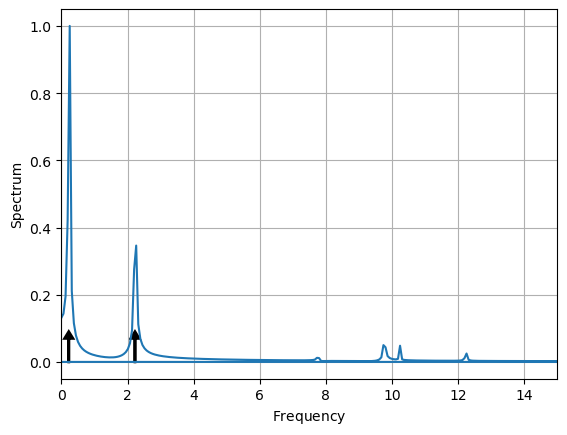}}
    \subfloat[Frequencies of $\bar{X}_1$]{\label{FFTXbar1_alp1eps01}\includegraphics[width=0.4\textwidth]{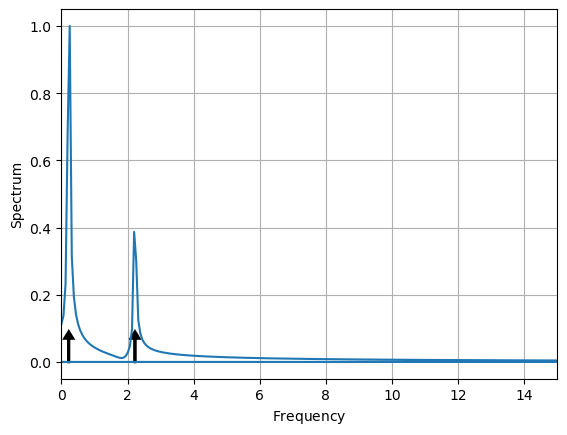}}
    \caption{\centering Explicit first order uniformly accurate scheme - There are two different frequencies in the sequence $\bar{X}_1$ (asymptotic model); the $\varepsilon-$model also has small frequencies due to perturbations arising from $\varepsilon$.}
    \label{fig:2freq}
\end{figure}

\begin{figure}
    \centering
    \subfloat[$(x_1,x_2)$ for $\varepsilon=0.1$]{\label{EP_x1x2_alp05eps01}\includegraphics[width=0.3\textwidth]{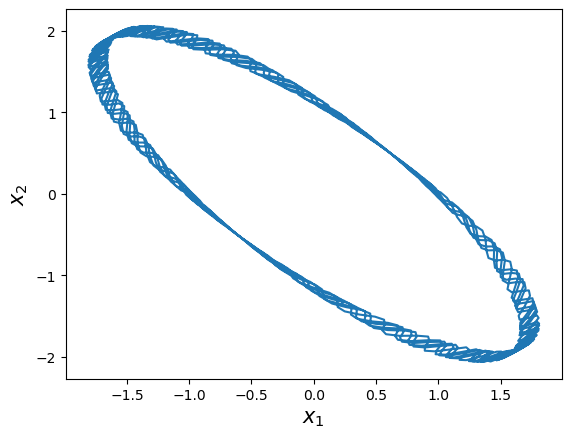}}
    \subfloat[$(x_1,x_2)$ for $\varepsilon=0.001$]{\label{EP_x1x2_alp05eps0001}\includegraphics[width=0.3\textwidth]{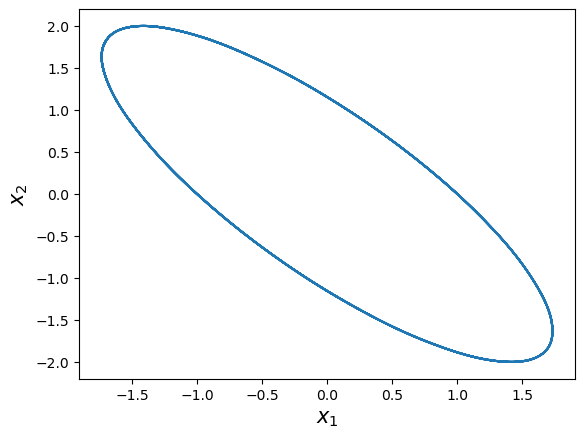}}
    \subfloat[Phase plot of $(\bar{X}_1,\bar{X}_2)$]{\label{EP_x1bx2b_alp05eps01}\includegraphics[width=0.3\textwidth]{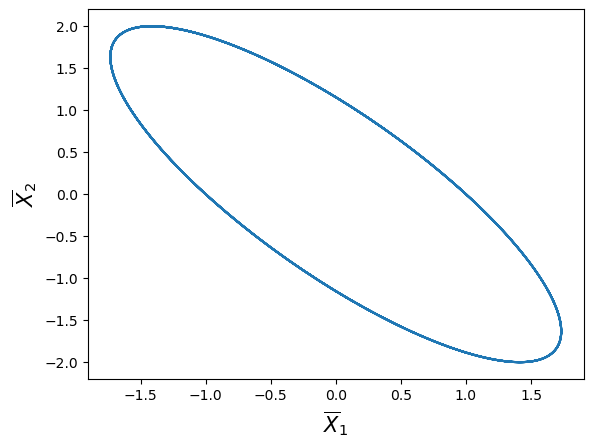}}

    \caption{\centering $B=0.5$, phase plots of $(x_1, x_2)$ for $\varepsilon=0.1, 0.001$ and of $(\bar{X}_1, \bar{X}_2)$ obtained with the energy preserving second order uniformly accurate scheme for $\dot{U}=A(t/\varepsilon)U$ with $A$ given by \eqref{A_numerical_section} and $\theta(s)=\cos(s)$. }
    \label{fig:1freq_alp05}
\end{figure}

\begin{figure}
    \centering
    \subfloat[$(x_1,x_2)$ for $\varepsilon=0.1$]{\label{EP_x1x2_alp1eps01}\includegraphics[width=0.3\textwidth]{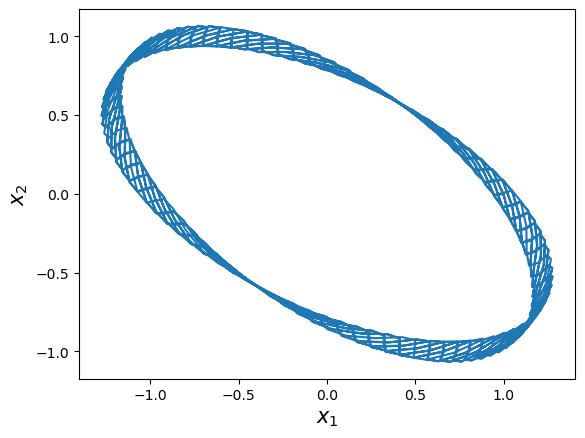}}
    \subfloat[$(x_1,x_2)$ for $\varepsilon=0.001$]{\label{EP_x1x2_alp1eps0001}\includegraphics[width=0.3\textwidth]{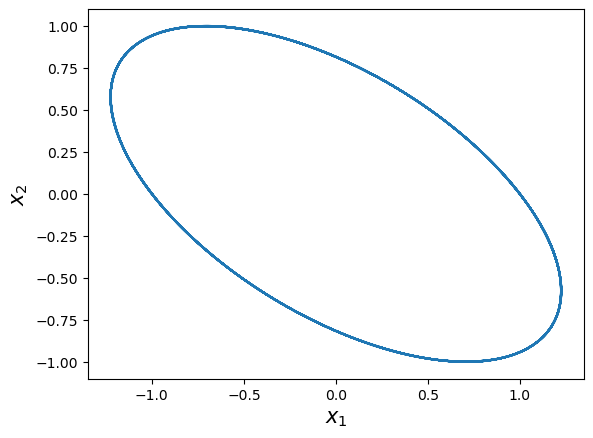}}
    \subfloat[Phase plot of $(\bar{X}_1,\bar{X}_2)$]{\label{EP_x1bx2b_alp1eps01}\includegraphics[width=0.3\textwidth]{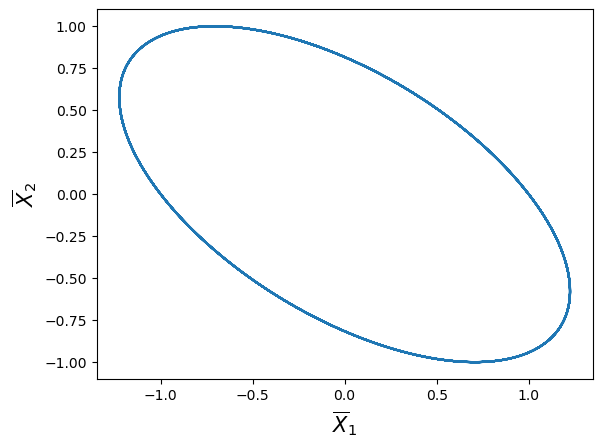}}

    \caption{\centering $B=1$, phase plots of $(x_1, x_2)$ for $\varepsilon=0.1, 0.001$ and of $(\bar{X}_1, \bar{X}_2)$ obtained with the energy preserving second order uniformly accurate scheme for $\dot{U}=A(t/\varepsilon)U$ with $A$ given by \eqref{A_numerical_section} and $\theta(s)=\cos(s)$. }
    \label{fig:1freq_alp1}
\end{figure}

\begin{figure}
    \centering
    \subfloat[$(x_1,x_2)$ for $\varepsilon=0.1$]{\label{EP_x1x2_alp5eps01}\includegraphics[width=0.3\textwidth]{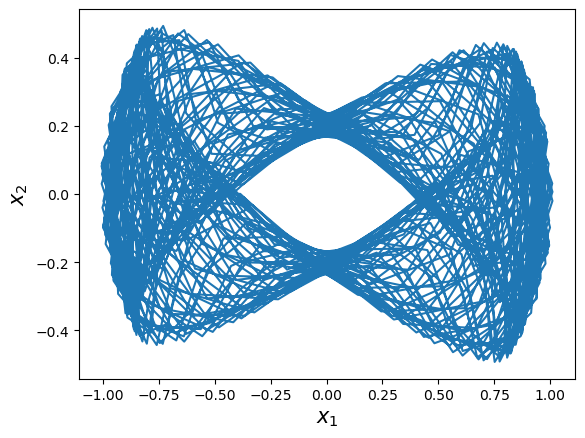}}
    \subfloat[$(x_1,x_2)$ for $\varepsilon=0.001$]{\label{EP_x1x2_alp5eps0001}\includegraphics[width=0.3\textwidth]{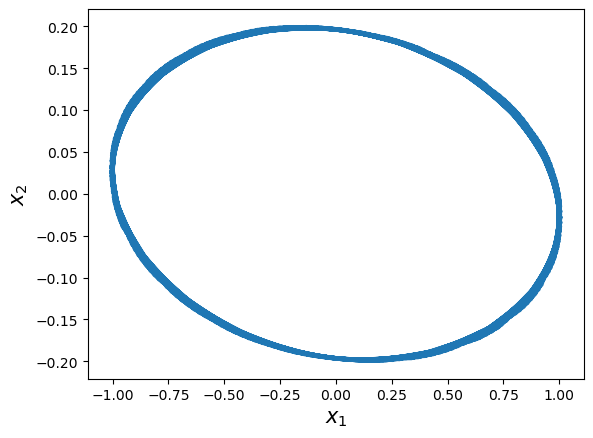}}
    \subfloat[Phase plot of $(\bar{X}_1,\bar{X}_2)$]{\label{EP_x1bx2b_alp5eps01}\includegraphics[width=0.3\textwidth]{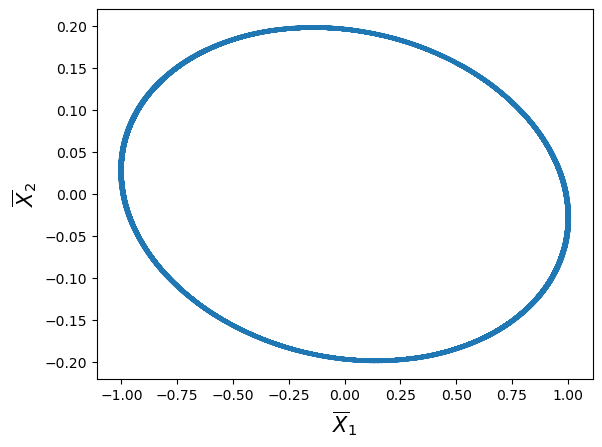}}

    \caption{\centering $B=5$, phase plots of $(x_1, x_2)$ for $\varepsilon=0.1, 0.001$ and of $(\bar{X}_1, \bar{X}_2)$ obtained with the energy preserving second order uniformly accurate scheme for $\dot{U}=A(t/\varepsilon)U$ with $A$ given by \eqref{A_numerical_section} and $\theta(s)=\cos(s)$.}
    \label{fig:1freq_alp5}
\end{figure}

Figure \ref{fig:2freq} corresponds to an explicit first order uniformly accurate scheme for the ODE system with $B=1,\Delta t = 0.001, \theta(s)=1+\cos{s}$, and the final time is $T=100$. We consider $\varepsilon=0.1$ and 
display the time history of $x_1$ together with its averaged counterpart 
$\bar{x}_1$ (see Figure \ref{fig:2freq}-(a)), the phase plot of $(x_1, x_2)$ (see Figure \ref{fig:2freq}-(b)) and the phase plot of $(\bar{x}_1, \bar{x}_2)$ (see Figure \ref{fig:2freq}-(c)). 
In Figure \ref{fig:2freq}-(d) (resp. (e)), we display the discrete Fourier transform 
of the sequence $(x_1)_n$ (resp. $(\bar{x}_1)_n$). Two main peaks 
can be observed corresponding to the 
main periods of the numerical solutions which are in very good agreement 
with the frequencies computed from the eigenvalues of $\langle A\rangle$ in the averaged model 
(equal to $i 0.224$ and $i 2.224$, corresponding to the arrows).
Let us remark that additional frequencies are observed for the spectrum of $(x_1)_n$ due to oscillations in $\varepsilon$. \\
\Cref{fig:1freq_alp05,fig:1freq_alp1,fig:1freq_alp5} show numerical results for the nonlinear problem obtained with the second order uniformly accurate energy preserving scheme with $\Delta t = 0.1, \theta(s)=\cos{s}$, for different values of $B=0.5,1$ and $5$ respectively. 
Different values of $\varepsilon=0.1,0.001$ are considered to study the behaviour of the system. Note that for this choice of $\theta$, there is only one frequency corresponding to the same imaginary value present in all the eigenvalues of $\langle A \rangle$ (equal to $i0.707107 B$ which is well captured by the numerical scheme, not shown). In 
\Cref{fig:1freq_alp05,fig:1freq_alp1,fig:1freq_alp5}, the phase plots of $(x_1, x_2)$ and $(\bar{x}_1, \bar{x}_2)$ are displayed to investigate the confinement property with respect to both $\varepsilon$ and $B$. 
We observe the domain of $(x_1,x_2)$ shrinks with increase in $B$ which illustrates an improved confinement. When $\varepsilon$ decreases, the perturbations in $(x_1,x_2)$ become smaller, resulting in better confinement. Further, energy is preserved up to machine accuracy (figures not shown) for the scheme \eqref{sav_nl_avg} of the averaged model.

\subsection{Vlasov-Poisson solver: Particle-In-Cell method}
\begin{figure}
    \centering
    \subfloat[$k_{1}=0.5$]{\label{alpx005kx05alpy0ky2pi_Eyn0}\includegraphics[width=0.3\textwidth]{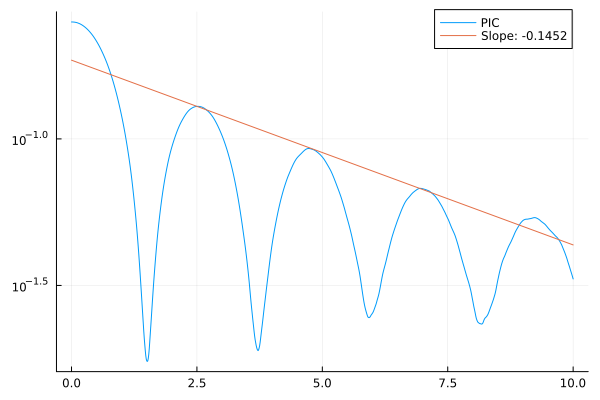}}
    \subfloat[$k_{1}=0.4$]{\label{alpx005kx04alpy0ky2pi_Eyn0}\includegraphics[width=0.3\textwidth]{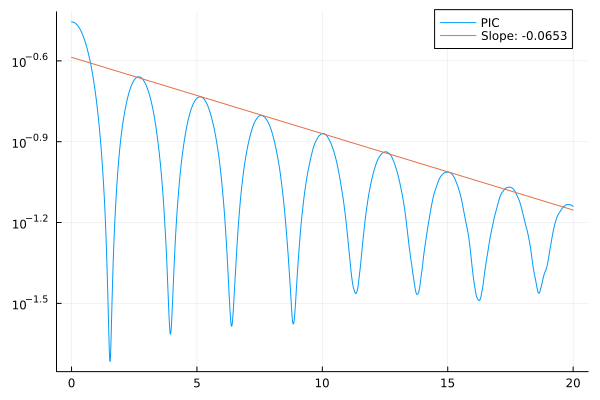}}
    \subfloat[$k_{1}=0.3$]{\label{alpx005kx03alpy0ky2pi_Eyn0}\includegraphics[width=0.3\textwidth]{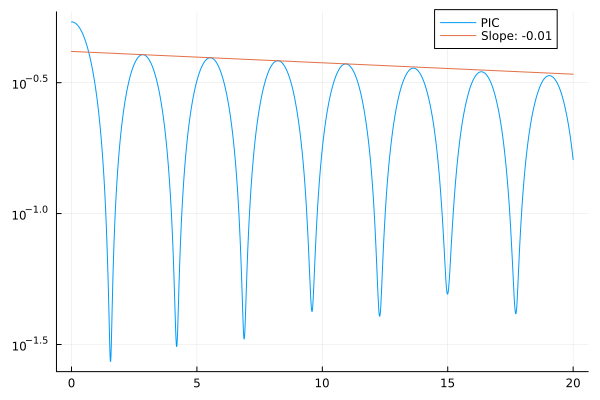}}
    \caption{\centering Vlasov-Poisson~\eqref{vlasov} with $B=0$. Landau damping test ($\xi_1=0.05, \xi_2=0$): time history of the electric energy for $k_1=0.5$ (left), $k_1=0.4$ (middle) and $k_1=0.3$ (right). }
    \label{fig:1DLDalp0}
\end{figure}

\begin{figure}
    \centering
    \subfloat[$B=0.01$]{\label{alpx005kx03alpy0ky2pi_Eyn0_lam001}\includegraphics[width=0.25\textwidth]{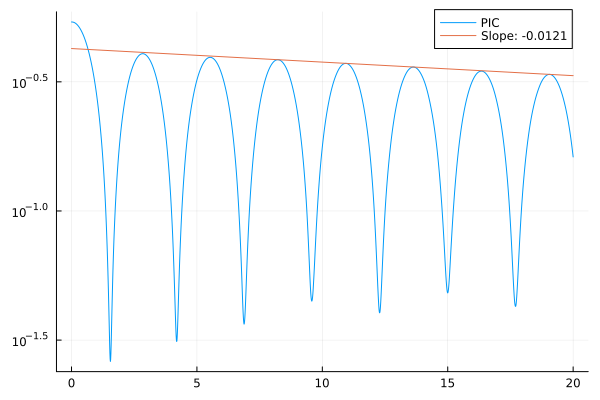}}
    \subfloat[$B=0.05$]{\label{alpx005kx03alpy0ky2pi_Eyn0_lam005}\includegraphics[width=0.25\textwidth]{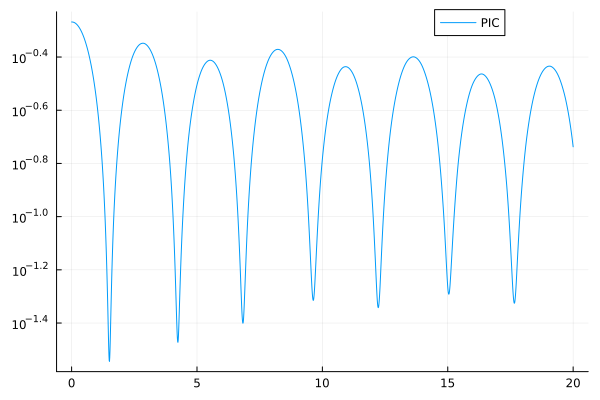}}
    \subfloat[$B=0.1$]{\label{alpx005kx03alpy0ky2pi_Eyn0_lam01}\includegraphics[width=0.25\textwidth]{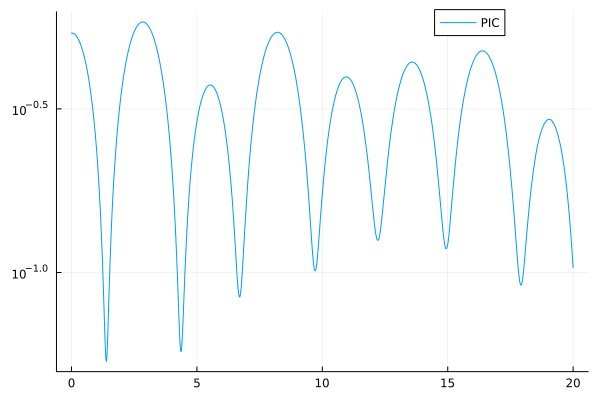}}
    \subfloat[$B=0.15$]{\label{alpx005kx03alpy0ky2pi_Eyn0_lam015}\includegraphics[width=0.25\textwidth]{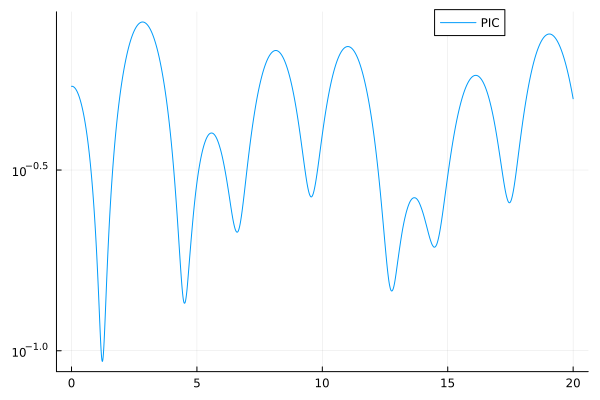}}
    \caption{\centering Vlasov-Poisson~\eqref{vlasov} with $\varepsilon=0.001$.  Landau damping test ($\xi_1=0.05, \xi_2=0, k_1=0.3$): time history of the electric energy for different values of $B$. }
    \label{fig:1DLDalpn0}
\end{figure}

\begin{figure}
    \centering
    \subfloat[$k_{1}=k_{2}=0.5$]{\label{alpx005kx05alpy005ky05_Eyn0}\includegraphics[width=0.3\textwidth]{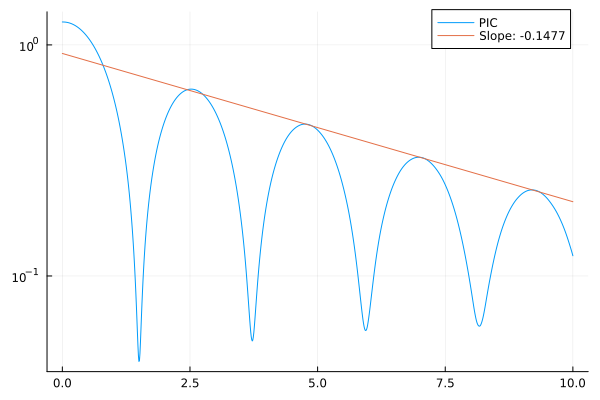}}
    \subfloat[$k_{1}=k_{2}=0.4$]{\label{alpx005kx04alpy005ky04_Eyn0}\includegraphics[width=0.3\textwidth]{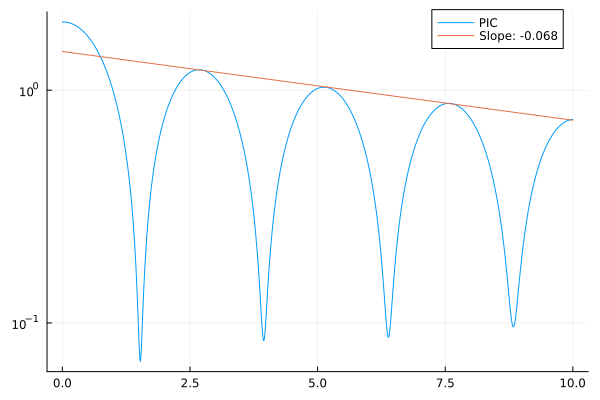}}
    \subfloat[$k_{1}=k_{2}=0.3$]{\label{alpx005kx03alpy005ky03_Eyn0}\includegraphics[width=0.3\textwidth]{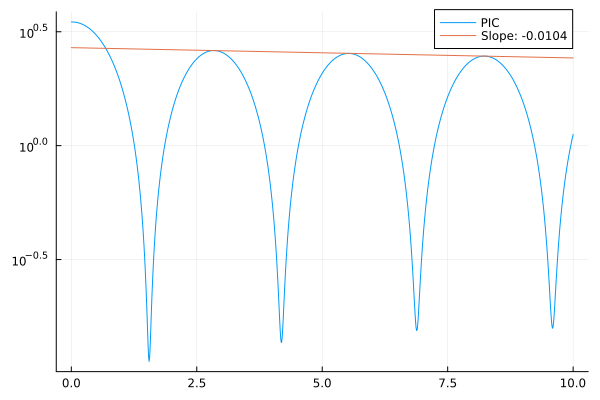}}
    \caption{\centering Vlasov-Poisson~\eqref{vlasov} with $B=0$.  Landau damping test ($\xi_1=\xi_2=0.05$): time history of the electric energy for $k_1=k_2=0.5$ (left), $k_1=k_2=0.4$ (middle) and $k_1=k_2=0.3$ (right). }
    \label{fig:2DLDalp0}
\end{figure}

\begin{figure}
    \centering
    \subfloat[$B=0.01$]{\label{alpx005kx03alpy005ky03_lam001}\includegraphics[width=0.25\textwidth]{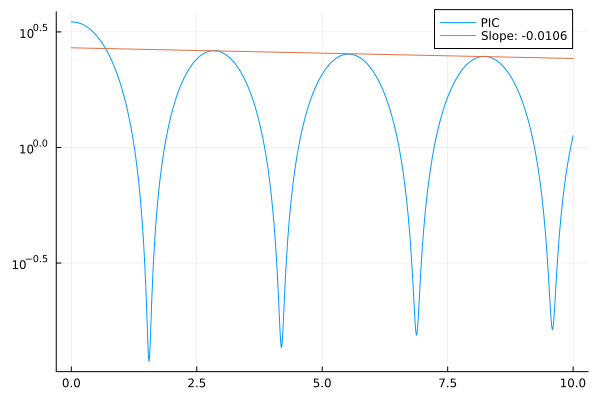}}
    \subfloat[$B=0.05$]{\label{alpx005kx03alpy005ky03_lam005}\includegraphics[width=0.25\textwidth]{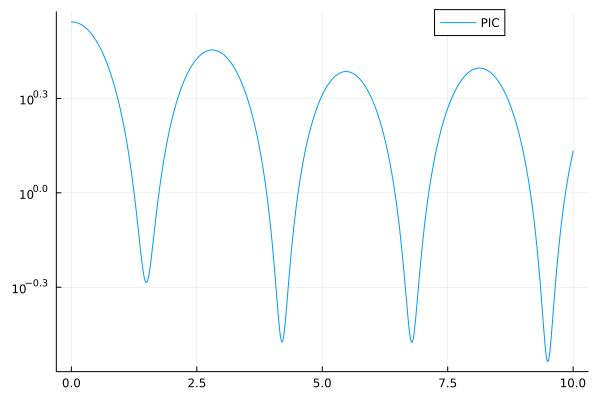}}
    \subfloat[$B=0.1$]{\label{alpx005kx03alpy005ky03_lam01}\includegraphics[width=0.25\textwidth]{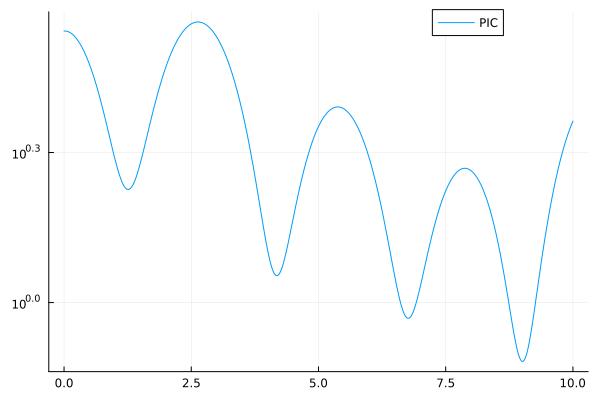}}
    \subfloat[$B=0.15$]{\label{alpx005kx03alpy005ky03_lam015}\includegraphics[width=0.25\textwidth]{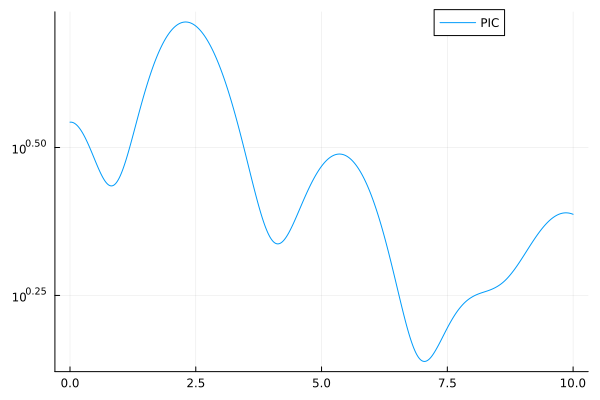}}
    \caption{\centering Vlasov-Poisson~\eqref{vlasov} with $\varepsilon=0.001$. Landau damping test ($\xi_1=\xi_2=0.05,  k_1=k_2=0.3$): time history of the electric energy for different values of $B$. }
    \label{fig:2DLDalpn0}
\end{figure}

In this section, we present some numerical results obtained for the Vlasov-Poisson equation by using PIC  method that employs our uniformly accurate numerical methods for solving the ODE system \eqref{charact_xq}. We consider the following initial condition 
\begin{equation}
\label{init_VP}
    f^{\text{in}}(x,v)= \Big(1+\xi_{1} \cos{(k_{1} x_1)}\Big) \Big(1+\xi_{2} \cos{(k_{2} x_2)}\Big) \frac{1}{2\pi} e^{-|v|^2/2}. 
\end{equation}
We also note that when $B=0$ in \eqref{charact_xq}, the problem degenerates to the Vlasov-Poisson equation. Thus, in this case, the initial condition \eqref{init_VP} corresponds to a Landau damping test which enables to validate our code by comparing our results with those in literature. Then, we investigate how the Landau damping behaves for non-zero values of $B$. We thus plot the time history electric energy $\int\int |E|^2dx_1 dx_2$ (semi-$\log$ scale).\\
In \Cref{fig:1DLDalp0,fig:1DLDalpn0}, we consider the one  dimensional case with $\xi_{2}=0$ and $k_2=2\pi$, whereas $\xi_{1}=0.05$. The domain is $x_1 \in[0,2\pi/k_{1}]$, $x_2 \in [0,2\pi/k_{2}]$. The number of  points along $x_1$ (resp. $x_2$) direction is $N_1=128$ (resp. $N_2=4$), 
 the number of particles is $N_p=100N_1N_2$, and time step is $\Delta t =0.01$. \Cref{fig:1DLDalp0} depicts the time history electric energy with $B=0$ for different values of $k_{1}$ 
 ($k_{1}=0.5,0.4,0.3$). The electric energy is expected to decrease exponentially fast in time with a rate which can be computed from the linear theory (see \cite{SonnendruckerBook}). In Figure \ref{fig:1DLDalp0}, it is observed that the obtained slopes are close to the theoretical values of damping observed in~\cite{SonnendruckerBook}. 
In \Cref{fig:1DLDalpn0}, we consider  $k_{1}=0.3$ case for non-zero values of $B=0.01,0.05,0.1$ and $0.15$, with magnetic field's oscillatory parameter $\varepsilon$ taken as $0.001$. When $B$ becomes larger $B\geq 0.1$, the damping disintegrates.  \\
In \Cref{fig:2DLDalp0,fig:2DLDalpn0}, we consider the two dimensional case  with $\xi_{1}=\xi_2=0.05$. The domain is $(x_1, x_2) \in[0,2\pi/k_{1}]\times [0,2\pi/k_{2}]$ (the numerical parameters are $N_1=N_2=128$, $N_p=100N_1N_2$ and $\Delta t=0.01$). \Cref{fig:2DLDalp0} depicts the Landau damping phenomena with $B=0$ for different values of $k_{1}=k_{2}$ such as $0.5,0.4,0.3$. There is not too much references for two-dimensional Landau damping test but we observe that the rate obtained with $k_1=k_2=0.5$ is very close to the theoretical one $-0.15$ (see~\cite{filbet_sonnen}). \\
In \cref{fig:2DLDalpn0}, we consider the $k_{1}=k_{2}=0.3$ case for non-zero values of $B=0.01,0.05,0.1$ and $0.15$, with magnetic field's oscillatory parameter $\varepsilon$ taken as $0.001$. Similar to the one dimensional case, the electric energy enjoys a very different behavior when $B$ increases.
 

\section{Conclusion}
In this work, we design and analyse uniformy accurate numerical schemes for highly oscillatory differential equations of the form $\dot{U}=A(t/\varepsilon)U + g(U)$. Both explicit and linearly implicit schemes are studied, according to the structure of the equation. In particular, when the averaged model enjoys a Hamiltonian structure, midpoint type schemes are proposed which degenerate into energy preserving schemes for the averaged model. Several numerical results illustrate the accuracy and geometric properties of the scheme. Further, the solvers have been used within PIC framework to solve the Vlasov-Poisson equation. \\
Several extensions to this work can be performed: for instance, the application to 
the dynamics of charged particles under the influence of strong magnetic fields, but also, links with Magnus integrators will be investigated. 

\appendix
\section{First order uniform accuracy of the standard midpoint scheme}
\label{appendix1}
In this first Appendix, we prove that the standard midpoint scheme~\eqref{mid_point_naive} is only first order UA. 
\begin{theorem}
The numerical scheme~\eqref{mid_point_naive} is first order uniformly accurate. 
\end{theorem}

\begin{proof}
Considering the error equation by subtracting~\eqref{duhamel_imex} with~\eqref{mid_point_naive}, we obtain
\begin{eqnarray}
    e_{n+1}-e_n &=& \int_{t_n}^{t_{n+1}} A(s/\varepsilon) \left(U(s)-\frac{U_n+U_{n+1}}{2} \right) ds \nonumber\\
    &=& \int_{t_n}^{t_{n+1}} A(s/\varepsilon) ds \left( \frac{e_n+e_{n+1}}{2} \right) + \int_{t_n}^{t_{n+1}} A(s/\varepsilon) \left(\frac{U(s)-U(t_n)}{2}+\frac{U(s)-U(t_{n+1})}{2} \right)  ds \nonumber\\
       \label{term_app}
    &\hspace{-3.5cm}=&\hspace{-2cm} \int_{t_n}^{t_{n+1}} \!\!\!\! A(s/\varepsilon) ds \left( \frac{e_n+e_{n+1}}{2} \right) 
    + \frac{1}{2} \int_{t_n}^{t_{n+1}} \!\!\!\! A(s/\varepsilon) \Big( \int_{t_n}^s A(s_1/\varepsilon) U(s_1) ds_1 \!-\! \int_s^{t_{n+1}} \!\!\!\! A(s_1/\varepsilon) U(s_1) ds_1 \Big) ds. 
\end{eqnarray}
Then, as the solution $U$ is uniformly bounded according to lemma~\ref{appendix_tech_lemma}-\ref{appendix_Gronwall}, the norm of $e_{n+1}$ becomes,
\begin{eqnarray*}
    \norm{e_{n+1}}\leq \frac{1+C\Delta t}{1-C\Delta t} \norm{e_n} +\frac{C}{1-C\Delta t} \Delta t^2.
\end{eqnarray*}
We then conclude by Gronwall lemma that $\norm{e_n}\leq C \Delta t$, with $C$ independent of $\Delta t$ and $\varepsilon$. 

However, it is well known that the midpoint scheme is second order for non stiff problems ($A$ independent of time or for $\varepsilon$ large). Let us thus focus on the last term in~\eqref{term_app}. It is possible to replace $U(s_1)$ by $U(t_n)+\int_{t_n}^{s_1} A(s_2/\varepsilon) U(s_2)ds_2$ to get 
\begin{eqnarray}
    e_{n+1}-e_n &=&\int_{t_n}^{t_{n+1}} \!\!\!\!A(s/\varepsilon) ds \left( \frac{e_n+e_{n+1}}{2} \right) + \frac{1}{2} \int_{t_n}^{t_{n+1}} \!\!\!\!A(s/\varepsilon) \Big( \int_{t_n}^s\!\!\! A(s_1/\varepsilon) ds_1  -\int_s^{t_{n+1}}\!\!\!\! A(s_1/\varepsilon)ds_1 \Big) ds \; U(t_n) \nonumber\\ 
    &&\hspace{-2cm}+ \frac{1}{2} \int_{t_n}^{t_{n+1}}\!\!\!\! A(s/\varepsilon) \Big( \int_{t_n}^s \!\! A(s_1/\varepsilon) 
     \int_{t_n}^{s_1} \!\!A(s_2/\varepsilon) U(s_2)ds_2) ds_1 \! -\! \int_s^{t_{n+1}} \!\!\!\! A(s_1/\varepsilon) \int_{t_n}^{s_1} \!\! A(s_2/\varepsilon) U(s_2)ds_2 ds_1 \Big) ds. \nonumber
\end{eqnarray}
The last term is of order $\Delta t^3$ uniformly in $\varepsilon$ due the uniform boundedness of $U$ and of $A$. The loss of uniform accuracy thus should come from the second term. A direct analysis gives a bound $\Delta t^2$ uniformly in $\varepsilon$, but this does not take the non stiff case into account.    
A Taylor expansion of $A(s_1/\varepsilon)$ gives 
$$
A(s_1/\varepsilon) = A(t_{n+1/2}/\varepsilon) + \frac{(s_1-t_{n+1/2})}{\varepsilon}A'(\xi/\varepsilon) 
$$
so that 
\begin{eqnarray*}
\int_{t_n}^{t_{n+1}}\!\!\!\! A(s/\varepsilon)\int_{t_n}^s A(s_1/\varepsilon)ds_1 ds&\!\!\!\!=\!\!\!\!& \int_{t_n}^{t_{n+1}}\!\!\!\! A(s/\varepsilon)(s-t_n)A(t_{n+1/2}/\varepsilon)ds +\int_{t_n}^{t_{n+1}}\!\!\!\! A(s/\varepsilon) \int_{t_n}^s \frac{(s_1-t_{n+1/2})}{\varepsilon}A'(\xi/\varepsilon)ds_1 ds\nonumber\\
\int_{t_n}^{t_{n+1}}\!\!\!\! A(s/\varepsilon)\int_s^{t_{n+1}} A(s_1/\varepsilon)ds_1 ds&\!\!\!\!=\!\!\!\!& \int_{t_n}^{t_{n+1}}\!\!\!\! A(s/\varepsilon)(t_{n+1}-s)A(t_{n+1/2}/\varepsilon)ds +\int_{t_n}^{t_{n+1}}\!\!\!\! A(s/\varepsilon) \int_s^{t_{n+1}} \frac{(s_1-t_{n+1/2})}{\varepsilon}A'(\xi/\varepsilon)ds_1 ds
\end{eqnarray*}
and the difference becomes 
\begin{eqnarray*}
\Big|\int_{t_n}^{t_{n+1}} A(s/\varepsilon) \Big[ \int_{t_n}^s A(s_1/\varepsilon) ds_1  -\int_s^{t_{n+1}} A(s_1/\varepsilon)ds_1 \Big] ds\Big| &\leq &  2\Big|\int_{t_n}^{t_{n+1}}\!\!\!\! A(s/\varepsilon)(s-t_{n+1/2})A(t_{n+1/2}/\varepsilon)ds\Big| \nonumber\\
&&\hspace{-6cm}+\Big|\int_{t_n}^{t_{n+1}}\!\!\!\! A(s/\varepsilon) \Big[\int_{t_n}^s \frac{(s_1-t_{n+1/2})}{\varepsilon}A'(\xi/\varepsilon)ds_1-\int_s^{t_{n+1}} \frac{(s_1-t_{n+1/2})}{\varepsilon}A'(\xi/\varepsilon)ds_1 \Big] ds\Big|.
\end{eqnarray*}
The last term is bounded by $C\Delta t^3/\varepsilon$. Let us now focus on the first term on the right hand side: 
\begin{eqnarray*}
\Big|\int_{t_n}^{t_{n+1}}\!\!\! A(s/\varepsilon)(s-t_{n+1/2})ds\Big|&=& \Big|\int_{t_n}^{t_{n+1}}\!\!\!\! \Big[A(t_{n+1/2}/\varepsilon) +\frac{(s-t_{n+1/2})}{\varepsilon}A'(\xi/\varepsilon)\Big] (s-t_{n+1/2})ds\Big| \nonumber\\ 
&=& \Big|\int_{t_n}^{t_{n+1}} \!\!\!\frac{(s-t_{n+1/2})^2}{\varepsilon}  A'(\xi/\varepsilon) ds\Big|  \leq  C\frac{\Delta t^3}{\varepsilon}. 
\end{eqnarray*}
We thus conclude with the following relation on the error 
$$
\|e_{n+1} \|\leq \frac{1+C\Delta t}{1-C\Delta t}\|e_n\| + \frac{C\Delta t^3/\varepsilon}{1-C\Delta t}, 
$$
from which we get, using Gronwall lemma $\|e_n\| \leq C\Delta t^2/\varepsilon$. 

Gathering the two estimates $\|e_n\| \leq C\Delta t$ and $\|e_n\| \leq C\Delta t^2/\varepsilon$ leads to $\|e_n\| \leq C\min(\Delta t,\Delta t^2/\varepsilon)$ and we conclude that the midpoint method is first order uniformly accurate. 
\end{proof}

\section{Technical lemmas}
\label{appendix_tech_lemma}
In this Appendix, we gather several technical lemmas which are used several times in different proofs. First, a very standard lemma  ensures that the solution of the ODE is uniformly bounded. 
\begin{lemma}[Uniform boundedness (with $\varepsilon$) of the exact solution]
\label{appendix_Gronwall}
Consider $U(t)$ the exact solution to~\eqref{lin_general}, where $g(t,U(t))$ is a Lipschitz continuous function with Lipschitz constant $K$ and $A$ is a bounded matrix-valued function. Then there exist positive constants $K_1$ and $K_2$, independent of $\varepsilon$, such that for $t\geq 0$, $\norm{U(t)} \leq K_1 \exp{\brac{\int_{0}^{t} K_2 ds }}$.
\end{lemma}
The proof is classical and derives from the Gronwall lemma and the following representation formula:
\begin{equation*}
    U(t)=U(0) + \int_{0}^{t} A\brac{s/\varepsilon} U(s) ds  + \int_{0}^{t} g(s,U(s)) ds.
\end{equation*}

The asymptotic expansions of averaged value for highly oscillating functions is useful. 
\begin{lemma}
\label{lemma_limit}
Let $\theta$ be a P-periodic smooth function and set $\langle \theta\rangle = \frac{1}{P} \int_{0}^P \theta(s) ds$ as its average. Then the following asymptotic expansions are satisfied ($t_n=n\Delta t$ for any integer $n$):
\begin{eqnarray}
\label{id1}
\int_a^b\theta(t/\varepsilon) dt &=& (b-a)\langle \theta\rangle + {\cal O}((b-a)\,\varepsilon), \\
\label{id2}
\int_{t_n}^{t_{n+1}} \int_{t_n}^t \theta(t_1/\varepsilon)dt_1 dt &=& \frac{\Delta t^2}{2}\langle \theta\rangle  +{\cal O}(\Delta t^2\,\varepsilon), \\
\label{id3}
\int_{t_n}^{t_{n+1}} \int_t^{t_{n+1}} \theta(t_1/\varepsilon)dt_1 dt &=& \frac{\Delta t^2}{2} \langle \theta\rangle +{\cal O}(\Delta t^2\,\varepsilon), \\
\label{id4}
\int_{t_n}^{t_{n+1}}\theta(t/\varepsilon) (t-t_n) dt &=& \frac{\Delta t^2}{2} \langle \theta\rangle +{\cal O}(\Delta t\,\varepsilon),\\ 
\label{id5}
\int_{t_n}^{t_{n+1}}\theta(t/\varepsilon) (t_{n+1}-t) dt &=&  \frac{\Delta t^2}{2} \langle \theta\rangle +{\cal O}(\Delta t\,\varepsilon). 
\end{eqnarray}
\end{lemma}
\begin{proof}
Without loss of generality, we may assume the period $P=2\pi$ and consider then the Fourier series expansion of $\theta(s)$, with coefficients $C_k =  \frac{1}{2\pi} \int_{0}^{2\pi} \theta(s) \mathsf{e}^{-iks} ds$, $k\in\mathbb{Z}$, having appropriate summability properties inherited from the smoothness of $\theta$ :
\begin{eqnarray*}
    \theta(s) = \sum_{k \in \mathbb{Z}} C_k \mathsf{e}^{iks} 
    = \langle \theta\rangle + \sum_{k \in \mathbb{Z}\setminus\{0\}} C_k \mathsf{e}^{ik s}, 
\end{eqnarray*}
we obtain
\begin{eqnarray*}
    \int_{a}^b \left( \theta(t/\varepsilon) - \langle \theta\rangle \right) dt = \sum_{k \in \mathbb{Z}\setminus\{0\}} C_k \int_{a}^b \mathsf{e}^{ikt/\varepsilon}   dt 
    = \varepsilon \sum_{k \in \mathbb{Z}\setminus\{0\}} \frac{1}{ik}C_k \left[ \mathsf{e}^{ikb/\varepsilon} - \mathsf{e}^{ika/\varepsilon} \right] = \mathcal{O}((b-a)\,\varepsilon).
\end{eqnarray*}
Thus, the identities~\eqref{id1},~\eqref{id2} and \eqref{id3} become evident:
\begin{gather*}
    \int_a^b\theta(t/\varepsilon) dt = (b-a)\langle \theta\rangle +{\cal O}((b-a)\,\varepsilon), \\
    \int_{t_n}^{t_{n+1}} \int_{t_n}^t \theta(t_1/\varepsilon)dt_1 dt = \int_{t_n}^{t_{n+1}} [\langle \theta\rangle (t-t_n) +{\cal O}(\Delta t\,\varepsilon)]dt  = \frac{\Delta t^2}{2}\langle \theta\rangle  +{\cal O}(\Delta t^2\,\varepsilon), \\
    \int_{t_n}^{t_{n+1}} \int_t^{t_{n+1}} \theta(t_1/\varepsilon)dt_1 dt = \int_{t_n}^{t_{n+1}} [\langle \theta\rangle (t_{n+1}-t)+{\cal O}(\Delta t\,\varepsilon)] dt  = \frac{\Delta t^2}{2} \langle \theta\rangle +{\cal O}(\Delta t^2\,\varepsilon).
\end{gather*}
To prove~\eqref{id4} and~\eqref{id5}, we use similar calculations in the following way
\begin{eqnarray*}
    \int_{t_n}^{t_{n+1}}(\theta(t/\varepsilon)-\langle \theta\rangle) (t-t_n) dt &=& \int_{t_n}^{t_{n+1}} \sum_{k \in \mathbb{Z}\setminus\{0\}} C_k \mathsf{e}^{ikt/\varepsilon} (t-t_n) dt \\
    \hspace{-7cm}&\hspace{-2cm}&\hspace{-2cm} =\varepsilon \sum_{k \in \mathbb{Z}\setminus\{0\}} \frac{1}{ik}C_k \left[ \mathrm{e}^{ik t/\varepsilon}(t-t_n) \right]_{t_n}^{t_{n+1}} - \varepsilon \int_{t_n}^{t_{n+1}} \sum_{k \in \mathbb{Z}\setminus\{0\}} \frac{1}{ik}C_k \mathrm{e}^{i k t/\varepsilon} dt = {\cal O}(\varepsilon \Delta t). 
\end{eqnarray*}
\end{proof}

\bibliographystyle{acm}
\bibliography{references}
\end{document}